%%%%%%%%%%%%%%%%%%%%%%%%%%%%%%%%%%%%%%%%%%%%%%%%%%%%%%%%%%%%%%%%%%%%%%%
%
% Artikel: Weak multiplier Hopf algebras II. The source and target algebras.
%
% Revised ersion 
%
%
%%%%%%%%%%%%%%%%%%%%%%%%%%%%%%%%%%%%%%%%%%%%%%%%%%%%%%%%%%%%%%%%%%%%%%%

\input amstex  % Uncomment indien nodig

\input amssym
\input amssym.def

\magnification 1100
\loadmsbm
\parindent 0 cm

\define\nl{\bigskip\item{}}
\define\snl{\smallskip\item{}}
\define\inspr #1{\parindent=20pt\bigskip\bf\item{#1}}
\define\iinspr #1{\parindent=27pt\bigskip\bf\item{#1}}
\define\einspr{\parindent=0cm\bigskip}

\define\tussen{\quad\qquad\qquad\quad}
\define\tussenen{\quad\qquad\text{and}\qquad\quad}

\define\ot{\otimes}

\define\tr{\triangleright}
\define\tl{\triangleleft}

\centerline{\bf Weak Multiplier Hopf Algebras II}
\snl\centerline{The source and target algebras}
\bigskip
\centerline{\it  Alfons Van Daele \rm $^{(1)}$ and \it Shuanhong Wang \rm $^{(2)}$}
\bigskip\bigskip
{\bf Abstract} 
\nl 
Let $(A,\Delta)$ be a {\it weak multiplier Hopf algebra} as introduced in [VD-W3] (see also [VD-W2]). It is a pair of a non-degenerate algebra $A$, with or without identity, and a coproduct $\Delta$ on $A$, satisfying certain properties. If the algebra has an identity and the coproduct is unital, then we have a Hopf algebra. If the algebra has no identity, but if the coproduct is non-degenerate (which is the equivalent of being unital if the algebra has an identity), then the pair would be a multiplier Hopf algebra. If the algebra has an identity, but the coproduct is not unital, we have a weak Hopf algebra. In the general case, we neither assume $A$ to have an identity nor do we assume $\Delta$ to be non-degenerate and so we work with a {\it genuine} weak multiplier Hopf algebra. It is called {\it regular} if its antipode is a bijective map from $A$ to itself.
\snl
In this paper, we {\it continue the study of weak multiplier Hopf algebras}. We recall the notions of the source and target maps $\varepsilon_s$ and $\varepsilon_t$, as well as of the source and target algebras. Then we investigate these objects further. Among other things, we show that the canonical idempotent $E$ (which is eventually $\Delta(1)$) belongs to the multiplier algebra $M(B\ot C)$ where $B=\varepsilon_s(A)$ and $C=\varepsilon_t(A)$ and that it is a {\it separability idempotent} (as studied in [VD4.v2]).  
If the weak multiplier Hopf algebra is regular, then also $E$ is a {\it regular} separability idempotent. 
\snl
We also consider {\it special cases and examples} in this paper. In particular, we see how for any weak multiplier Hopf algebra $(A,\Delta)$, it is possible to make $C\ot B$ (with $B$ and $C$ as above) into a new weak multiplier Hopf algebra. In a sense, it forgets the 'Hopf algebra part' of the original weak multiplier Hopf algebra and only remembers the source and target algebras. It is in turn generalized to the case of any pair of non-degenerate algebras $B$ and $C$ with  a separability idempotent $E\in M(B\ot C)$. We get another example using this theory associated to any discrete quantum group (a multiplier Hopf algebra of discrete type with a normalized cointegral).  Finally we also consider the well-known 'quantization' of the groupoid that comes from an action of a group on a set. All these constructions provide interesting new examples of weak multiplier Hopf algebras (that are not weak Hopf algebras). 
\nl 
Date: {\it 25 September 2015}. Revised version.
\vskip 1 cm
\hrule
\bigskip
\parindent 0.7 cm
\item{($1$)} Department of Mathematics, University of Leuven, Celestijnenlaan 200B, \newline
B-3001 Heverlee, Belgium. {\it E-mail}: Alfons.VanDaele\@wis.kuleuven.be
\item{($2$)} Department of Mathematics, Southeast University, Nanjing 210096, China. \newline
{\it E-mail}:  Shuanhwang2002\@yahoo.com {\it or} Shuanhwang\@seu.edu.cn 
\parindent 0 cm

\newpage

\bf 0. Introduction \rm
\nl
Consider a groupoid $G$. It is a set with a product that is not defined for all pairs of elements $p,q\in G$. Only if the so-called source $s(p)$ of $p$ equals the target $t(q)$ of $q$, then $pq$ is defined in $G$. The source and target are maps from $G$ to the so-called units of $G$. Here we consider the units as a subset of $G$.  The product is associative in the obvious sense and for any element $p\in G$, there is a unique inverse $p^{-1}$ characterized by the property that $s(p^{-1})=t(p)$ and $t(p^{-1})=s(p)$ and that $p^{-1}p=s(p)$ and $pp^{-1}=t(p)$.
\snl
We refer to basic works on groupoids for a more precise definition and details. See further in this introduction under the item {\it Basic references}.
\nl
With any  groupoid $G$ are associated two (regular) weak multiplier Hopf algebras (in duality).
\snl
First there is the algebra $A$, defined as the space $K(G)$ of complex functions on $G$ with finite support and pointwise product. A coproduct $\Delta$ on $K(G)$ is defined by  
$$\Delta(f)(p,q)=	
\cases
  f(pq) & \text{if } pq \text{ is defined},\\
  0 & \text{otherwise}.
\endcases
$$
The pair $(A,\Delta)$ is a regular weak multiplier Hopf algebra (in the sense of [VD-W3]). The idempotent multiplier $E$ in $M(A\ot A)$ (playing the role of $\Delta(1)$ and sometimes called the {\it canonical idempotent} of the weak multiplier Hopf algebra) is given by the function on pairs $(p,q)$ in $G\times G$ that is $1$ if $pq$ is defined and $0$ if this is not the case. The antipode $S$ is defined by $(S(f))(p)=f(p^{-1})$ whenever $f\in K(G)$ and $p\in G$.
\snl
In this example, the {\it source map} $\varepsilon_s$ from $A$ to $M(A)$ is defined by $(\varepsilon_s(f))(p)=f(p^{-1}p)$ whenever $p\in G$ and $f\in K(G)$. The multiplier algebra of the source algebra $\varepsilon_s(A)$ can be identified with the algebra $A_s$ of all functions on $G$ so that $f(p)=f(q)$ whenever $p,q\in G$ satisfy $s(p)=s(q)$. Similarly, the {\it target map} $\varepsilon_t$ from $A$ to $M(A)$ is defined by $(\varepsilon_t(f))(p)=f(pp^{-1})$ for all $p$ and $f\in K(G)$. The multiplier algebra of the target algebra $\varepsilon_t(A)$ will be the algebra $A_t$ of functions $f$ on $G$ so that $f(p)=f(q)$ if $t(p)=t(q)$ for $p,q\in G$.
Recall that the ranges $\varepsilon_s(A)$ and $\varepsilon_t(A)$ of the source and target maps can be strictly smaller than the algebras $A_s$ and $A_t$ respectively. This happens when the set of units is infinite.
\snl
We refer to Section 1 in [VD-W3] for more details on this example. See also Example 3.1 in Section 3 of this paper.
\nl
For the second case, we take the algebra $B$, defined as the groupoid algebra $\Bbb C G$ of $G$. If we use $p\mapsto \lambda_p$ for the canonical embedding of $G$ in $\Bbb C G$, then if $p,q\in G$, we have $\lambda_p\lambda_q=\lambda_{pq}$  if $pq$ is defined and $0$ otherwise. Here the canonical idempotent $E$ is given by $\sum \lambda_e\ot \lambda_e$ where the sum is only taken over the units $e$ of $G$. The antipode is given by $S(\lambda_p)=\lambda_{p^{-1}}$ for all $p\in G$.
\snl
The source and target maps are given by $\varepsilon_s(\lambda_p)=\lambda_e$ where $e=s(p)$ and $\varepsilon_t(\lambda_p)=\lambda_e$ where now $e=t(p)$. 
Here the multiplier algebras $B_s$ and $B_t$ of the source and target algebras $\varepsilon_s(B)$ and $\varepsilon_t(B)$ coincide and it is the multiplier algebra of the span of elements of the form $\lambda_e$ where $e$ is a unit of $G$.
Also in this case, the images of the source and target maps can be strictly smaller than the  algebras $B_s$ and $B_t$.
\snl
Again, we refer to Section 1 in [VD-W3] for more details on this example (and again also Example 3.1 in Section 3).
\nl
These two cases are dual to each other. The duality is given by $\langle f,\lambda_p\rangle=f(p)$ whenever $f\in K(G)$ and $p\in G$. We give more details (about this duality) in [VD-W5] where we treat duality for regular weak multiplier Hopf algebras with integrals.
\nl
In this paper we continue the study of the source and target maps $\varepsilon_s$ and $\varepsilon_t$ and their images, the source and target algebras $\varepsilon_s(A)$ and $\varepsilon_t(A)$,  for a general (possibly non-regular) weak multiplier Hopf algebra $(A,\Delta)$. The multiplier algebras of the source and target algebras can be embedded in the multiplier algebra $M(A)$ of $A$. In the case of a regular weak multiplier Hopf algebra, they have a nice concrete characterization. It is not clear if this is still possible in the non-regular case. 
\snl
Also in this paper, we construct certain examples of weak multiplier Hopf algebras. The constructions are known in the case of finite-dimensional weak Hopf algebras and it turns out that they can be formulated also in this more general framework. Of course, some care is needed because coproducts do not map into the tensor product $A\ot A$, but into the multiplier $M(A\ot A)$ of this tensor product. 
\snl
In this sense, the (sub)title of this paper is somewhat misleading and too restrictive. In earlier papers on the subject, we hardly looked at other than the trivial motivating examples coming from a groupoid or weak Hopf algebras. In this paper we take advantage of the further development of the theory to consider examples that are closely related and using the results we obtain on the source and target algebras. In our forthcoming paper on the subject [VD-W5], where we treat integrals and duality, we will use these examples again and apply duality to get still other examples of weak multiplier Hopf algebras. Because we do not have to restrict to the finite-dimensional case, we get many more interesting examples that do not fit into the original theory of weak Hopf algebras. 
\nl
\it Content of the paper \rm
\nl
In {\it Section} 1 we recall some of the basic notions and results on weak multiplier Hopf algebras as studied in our first papers on the subject ([VD-W2] and [VD-W3]). In particular, we will explain some of the covering properties as this will be important for the rest of the paper.
\snl
In the earlier papers on the subject, we briefly looked already at the source and target maps $\varepsilon_s$ and $\varepsilon_t$ and their images, the source and target algebras.  In {\it Section} 2 we investigate these objects further. We recall the definitions and some of the basic properties that are found already in [VD-W3]. Notice that we make a {\it change in terminology}. We will now call the image  $\varepsilon_s(A)$ of the source map the source algebra and the image $\varepsilon_t(A)$ of the target map the target algebra. In [VD-W4.v1] we used these terms for the multiplier algebras that can be characterized nicely in the regular case. Because now we are also studying the non-regular case, these multiplier algebras no longer seem to have the same characterization and this is what motivated us to change this terminology. We comment more on this in Section 2.
\snl
Indeed, in the regular case, we show that the multiplier algebras $M(\varepsilon_s(A))$ and $M(\varepsilon_t(A))$ of the images $\varepsilon_s(A)$ and $\varepsilon_t(A)$ of the source and target maps can be nicely characterized as certain subalgebras of the multiplier algebra $M(A)$. 
\snl
In the general case, we show that the canonical idempotent $E$ has all the properties of a separability idempotent (as studied in [VD4.v2]). It turns out to be a regular one if the weak Hopf algebra is regular. Finally we use the various results to show that the underlying algebra $A$ of any weak multiplier Hopf algebra $(A,\Delta)$ has local units. Recall that in [VD-W3], we only could show this in the regular case. 
\snl
In {\it Section} 3 we study special cases and examples. We start again with the two examples associated with a groupoid. We will be very short here as we include this mainly for completeness. These examples have been considered in earlier papers (see e.g.\ [VD-W3]). Then we consider any 
weak multiplier Hopf algebra $(A,\Delta)$ and we associate a new weak multiplier Hopf algebra $(P,\Delta_P)$ where the underlying algebra $P$ is $\varepsilon_t(A)\ot \varepsilon_s(A)$ and the coproduct is given by the formula 
$$\Delta_P(c\ot b)=c\ot E\ot b$$ 
for $b\in \varepsilon_s(A)$ and $c\in \varepsilon_t(A)$ and where $E$ is the canonical multiplier in $M(A\ot A)$. We also use this example further as a model for the construction of an abstract version of this case. Then we take any pair of non-degenerate algebras $B$ and $C$ and start with a so-called separability idempotent $E$ in the multiplier algebra $M(B\ot C)$. We take $P=C\ot B$ and $\Delta_P$ as above. These two examples are 'quantizations' of the trivial groupoid $G$ constructed from a set $X$ by taking $G=X\times X$ with product $(z,y)(y,x)=(z,x)$ when $x,y,z\in X$. 
\snl
This groupoid in turn is related with the case of a groupoid $G$ constructed from a (left) action of a group $H$ on a set $X$. Now $G$ consists of triples $(y,h,x)$ where $x,y\in X$ and $h\in H$ and $y=h\tr x$ and where $\tr$ is used to denote the action. The product is given by $(z,k,y)(y,h,x)=(z,kh,x)$. And finally, also this groupoid will be quantized (at least in a certain sense to be explained in this section).
\snl
The starting point is again a pair of non-degenerate algebras $B$ and $C$ with a separability idempotent $E$ in the multiplier algebra $M(B\ot C)$. Moreover there is a (regular) multiplier Hopf algebra $Q$ that acts from the right on $B$ and from the left on $C$ in such a way that $B$ is a right $Q$-module algebra and $C$ a left $Q$-module algebra. These objects are related with the requirement that the right action of $Q$ on $C$ induces via $E$ the left action of $Q$ on $B$. See Section 3 for a more precise statement. The two-sided smash product $P$ is defined as the algebra generated by $B$, $C$ and $Q$ with $B$ and $C$ commuting and the commutation rules between $B$ and $Q$ determined by the left action of $Q$ on $B$ and the ones between $C$ and $Q$ determined by the right action of $Q$ on $C$. It carries a natural coproduct making $P$ into a weak multiplier Hopf algebra. 
\snl
Finally, in {\it Section} 4 we draw some conclusions and discuss possible further research on this subject.
\nl
In a forthcoming paper on the subject, we treat integrals and duality  and we consider more examples (cf.\ [VD-W5]). 
\snl
The material studied in this paper is closely related with the theory of {\it (regular) multiplier Hopf algebroids}, as developed in [T-VD1], where the theory of weak multiplier Hopf algebras is treated within an algebroid framework. See also [T-VD2] for the relation between the two concepts.
\snl
We  also like to refer to the paper on {\it weak multiplier bialgebras} by B\"ohm, G\'omez-Torecillas and L\'opez-Centella (see [B-G-L]) where the notion of a weak multiplier bialgebra is developed. In this theory, the source and target maps, as well as the source and target algebras, play a crucial role. See also [K-VD] where a Larson-Sweedler type theorem is proven for these weak multiplier bialgebras.
\nl
\it Conventions and notations \rm
\nl
We only work with algebras $A$ over $\Bbb C$ (although we believe that this is not essential and that it is possible to obtain the same results for algebras over other, more general fields). We do not assume that they are unital but we need that the product is non-degenerate. We also assume our algebras to be idempotent (that is $A^2=A$). In fact, it turns out that the algebras we encounter in this theory always have local units. We have seen this already in [VD-W3], in the regular case. Then of course, the product is automatically non-degenerate and also the algebra is idempotent. 
\snl
When $A$ is such an algebra, we use $M(A)$ for the multiplier algebra of $A$. When $m$ is in $M(A)$, then by definition we can define $am$ and $mb$ in $A$ for all $a,b\in A$ and we have $(am)b=a(mb)$. The algebra $A$ sits in $M(A)$ as an essential two-sided ideal and $M(A)$ is the largest algebra with identity having this property. 
\snl
Recall that a homomorphism $\gamma:A\to M(B)$, where $A$ and $B$ are non-degenerate algebras, is called non-degenerate if $\gamma(A)B=B$ and $B\gamma(A)=B$. In that case, there is a unique extension of $\gamma$, still denoted by $\gamma$, to a unital homomorphism from $M(A)$ to $M(B)$. There is a similar result for non-degenerate anti-homomorphisms.
\snl
We consider $A\ot A$, the tensor product of $A$ with itself. It is again an idempotent, non-degenerate algebra and we can consider the multiplier algebra $M(A\ot A)$. The same is true for a multiple tensor product. We use $\zeta$ for the flip map on $A\ot A$, as well as for its natural extension to $M(A\ot A)$.
\snl
We use $1$ for the identity in any of these  multiplier algebras. On the other hand, we mostly use $\iota$ for the identity map on $A$ (or other spaces), although sometimes, we also write $1$ for  this map. The identity element in a group is denoted by $e$. If $G$ is a groupoid, we will also use $e$ for units. Units are considered as being elements of the groupoid and we use $s$ and $t$ for the source and target maps from $G$ to the set of units. 
\snl
When $A$ is an algebra, we denote by $A^{\text{op}}$ the algebra obtained from $A$ by reversing the product. When $\Delta$ is a coproduct on $A$, we denote by $\Delta^{\text{cop}}$ the coproduct on $A$ obtained by composing $\Delta$ with the flip map $\zeta$.
\snl
For a coproduct $\Delta$, as we define it in Definition 1.1 of [VD-W3], we assume that $\Delta(a)(1\ot b)$ and $(a\ot 1)\Delta(b)$ are in $A\ot A$ for all $a,b\in A$. This allows us to make use of the {\it Sweedler notation} for the coproduct. The Sweedler notation is first explained in [Dr-VD], but only for the case of regular coproducts. In [VD3] an approach is developed in the case where the underlying algebras have local units. In the more recent paper [VD6], this condition is not assumed. However, it should be mentioned that the Sweedler notation is essentially just what is says, a notation. It is a way to denote formulas in a more transparent way. This point of view is explained in [VD6] and the reader is advised to look at that note for understanding the use of the Sweedler notation for weak multiplier Hopf algebras as in this paper.
\nl
\it Basic references \rm
\nl
For the theory of Hopf algebras, we refer to the standard works of Abe [A] and Sweedler [S]. For multiplier Hopf algebras and integrals on multiplier Hopf algebras, we refer to [VD1] and [VD2]. Weak  Hopf algebras have been studied in [B-N-S] and [B-S] and more results are found in [N] and [N-V1]. Various other references on the subject can be found in [Va]. In particular, we refer to [N-V2] because we will use notations and conventions from this paper when dealing with weak Hopf algebras.
\snl
For the theory of groupoids, we refer to [Br], [H], [P] and [R]. 
\nl
\bf Acknowledgments \rm
\nl
The first named author (Alfons Van Daele) would like to express his thanks to his coauthor Shuanhong Wang for motivating him to start the research on weak multiplier Hopf algebras and for the hospitality when visiting the University of Nanjing in 2008 for the first time and later again in 2012 and 2014. 
\snl
The second named author (Shuanhong Wang) would like to thank his coauthor for his help and advice when visiting the Department of Mathematics of the University of Leuven in Belgium during several times the past years. This work is partially supported by the NSF of China (No 11371088) and the NSF of Jiangsu Province, China (No. BK 2012736).
\snl
We also like to thank an anonymous referee who pointed out some shortcomings of an earlier version of this paper (cf.\ reference [VD-W4.v1]). This motivated us to do some more research and it resulted in an improvement of the results, mainly for the non-regular case.
\nl\nl

\bf 1. Preliminaries on weak multiplier Hopf algebras \rm
\nl
Let $(A,\Delta)$ be {\it a weak multiplier Hopf algebra} as in Definition 1.14 of [VD-W3]. In general, we do not assume that it is regular. On the other hand, we also recall some of the results that are only true in the regular case.
\snl
$A$ is {\it an algebra} over $\Bbb C$, with or without identity but with a product that is non-degenerate (as a bilinear map). The algebra is also idempotent in the sense that $A=A^2$ (meaning that any element in $A$ is a sum of products of elements of $A$). In Proposition 4.9 of [VD-W3], we showed that in the regular case, the underlying algebra automatically has local units. In fact, the result turns out to be true also in the non-regular case. We will obtain a proof in this paper (see Proposition 2.21 in Section 2). Remark that for an algebra with local units, the product is automatically non-degenerate and the algebra is idempotent.
\nl
There is a {\it coproduct} $\Delta$ on $A$. It is a homomorphism from $A$ to the multiplier algebra $M(A\ot A)$ of the tensor product $A\ot A$ of $A$ with itself. It is not assumed that it is non-degenerate (see further). The {\it canonical maps} $T_1$, $T_2$, $T_3$ and $T_4$ are  linear maps defined on $A\ot A$ by
$$\align 
T_1(a\ot b)&=\Delta(a)(1\ot b)
\qquad\quad\text{}\qquad\quad
T_2(c\ot a)=(c\ot 1)\Delta(a)\\
T_3(a\ot b)&=(1\ot b)\Delta(a)
\qquad\quad\text{}\qquad\quad
T_4(c\ot a)=\Delta(a)(c\ot 1).
\endalign$$
In general, it is assumed that $T_1$ and $T_2$ have range in $A\ot A$. If also  $T_3$ and $T_4$ map into $A\ot A$, then the {\it coproduct} is called {\it regular}.
\snl
The coproduct is assumed to be {\it full}. This means that the smallest subspaces $V$ and $W$ of $A$ satisfying
$$\Delta(a)(1\ot b)\in V\ot A
\qquad\quad\text{and}\qquad\quad
(c\ot 1)\Delta(a)\in A\ot W$$
for all $a,b\in A$ are $A$ itself. If the coproduct is regular, then a similar property will also be true for the maps $T_3$ and $T_4$ and so both the flipped coproduct $\Delta^{\text{cop}}$ on $A$ and the original coproduct on $A^{\text{op}}$ will also be full coproducts. 
\snl
Fullness of the coproduct implies that any element in $A$ is a linear span of elements of the form $(\iota\ot\omega)(\Delta(a)(1\ot b))$ where $a,b\in A$ and where $\omega$ is a linear functional on $A$. Similarly for the span of elements $(\omega\ot\iota)((c\ot 1)\Delta(a))$ with $a,c\in A$ and a linear functional $\omega$ on $A$. In fact, this property is equivalent with fullness of the coproduct. We have a result of the same type for fullness of a regular coproduct. See e.g.\ Proposition 1.6 in [VD-W1] and also Lemma 1.11 in [VD-W2]
\snl
Furthermore, it is assumed that there is {\it a counit}. This is a linear map $\varepsilon:A\to \Bbb C$ satisfying
$$(\varepsilon\ot\iota)(\Delta(a)(1\ot b))=ab
\qquad\quad\text{and}\qquad\quad
(\iota\ot\varepsilon)((c\ot 1)\Delta(a))=ca$$
for all $a,b,c$ in $A$. 
Similar formulas will be true for the other canonical maps in the case of a regular coproduct. 
\snl
Because the coproduct is assumed to be full, this counit is unique in the following sense. 
Assume that $\varepsilon$ and $\varepsilon'$ are linear maps such that 
$$(\varepsilon\ot\iota)(\Delta(a)(1\ot b))=ab
\qquad\quad\text{and}\qquad\quad
(\iota\ot\varepsilon')((c\ot 1)\Delta(a))=ca$$
for all $a,b,c$ in $A$. Then already $\varepsilon=\varepsilon'$. This is proven by applying $\iota\ot\varepsilon\ot\iota$ on the right hand side  and $\iota\ot\varepsilon'\ot\iota$ on the left hand side of the equation that expresses coassociativity of the coproduct
$$(c\ot 1\ot 1)(\Delta\ot\iota)(\Delta(a)(1\ot b))=(\iota\ot\Delta)((c\ot 1)\Delta(a))(1\ot 1\ot b).$$
In the two cases we get the same result, namely 
$(c\ot 1)\Delta(a)(1\ot b)$.
This is true for all $a,b,c\in A$ and from the fullness of the coproduct, it follows that $\varepsilon=\varepsilon'$.
\snl
It is not clear if there is a uniqueness result without the assumption that the coproduct is full. And it is also not clear if the existence of a counit, in the non-unital case, implies fullness of the coproduct. Remark that in general, the counit is not a homomorphism in the case of weak multiplier Hopf algebras. 
\snl
It seems not possible to construct a counit, even given that the coproduct is full. Therefore, the existence of the counit is part of the axioms for weak multiplier Hopf algebras.
\nl
There is an idempotent element $E$ in $M(A\ot A)$, called the {\it canonical idempotent}, giving the ranges of the canonical maps $T_1$ and $T_2$ as 
$$ \Delta(A)(1\ot A)=E(A\ot A)
\qquad\quad\text{and}\qquad\quad
(A\ot 1)\Delta(A)=(A\ot A)E.
$$
If the weak multiplier Hopf algebra is regular, we also have these properties for the ranges of the canonical maps $T_3$ and $T_4$. So in that case, we also have
$$ \Delta(A)(A\ot 1)=E(A\ot A)
\qquad\quad\text{and}\qquad\quad
(1\ot A)\Delta(A)=(A\ot A)E
$$
with the same idempotent. This element is uniquely determined and it satisfies 
$$\Delta(a)E=\Delta(a)=E\Delta(a)$$ 
for all $a\in A$.
\snl
We see that the coproduct is degenerate if $E$ is strictly smaller than $1$. However, still the coproduct can be extended in a unique way to a homomorphism from $M(A)$ to $M(A\ot A)$ (again denoted by $\Delta$) provided we assume $\Delta(1)=E$. Similarly, the homomorphisms $\Delta\ot\iota$ and $\iota\ot\Delta$ have unique extension to $M(A\ot A)$ such that, again using the same symbols for these extensions, we have 
$$(\Delta\ot\iota)(1)=E\ot 1
\qquad\quad\text{and}\qquad\quad
(\iota\ot\Delta)(1)=1\ot E.$$
We use $1$ for the identity, both in $M(A)$ and in $M(A\ot A)$. We have  $(\Delta\ot\iota)(E)= (\iota\ot\Delta)(E)$. It is further assumed that 
$$(\Delta\ot\iota)(E)=(E\ot 1)(1\ot E)=(1\ot E)(E\ot 1).$$
The last equality means, in a sense that can be made precise, that the left and the right legs of $E$ commute.
\snl
The left and the right legs of $E$ are also {\it big enough} in the following sense. 

\inspr{1.1} Lemma \rm
If $a\in A$ and if $E(1\ot a)=0$, then $a=0$. Similarly $a=0$, if $(1\ot a)E=0$, $(a\ot 1)E=0$ or if $E((a\ot 1)=0$.

\bf\snl Proof\rm:
Assume $a\in A$. If $E(1\ot a)=0$ then $\Delta(b)(1\ot a)=0$ for all $b\in A$. If we apply the counit $\varepsilon$ on the first leg of this equality we find $ba=0$ for all $b$ and so $a=0$. If $E(a\ot 1)=0$ we get $(c\ot 1)\Delta(b)(a\ot 1)=0$ for all $b,c\in A$. Now we apply the counit on the second leg and we find $cba=0$ for all $b,c\in A$. Again this implies $a=0$. A similar argument works for the two other cases.
\hfill $\square$ \einspr

There is a {\it unique antipode} $S$. It is a linear map from $A$ to the multiplier algebra $M(A)$. It is  an anti-algebra map in the sense that $S(ab)=S(b)S(a)$ for all $a,b\in A$ and it is an anti-coalgebra map meaning that $\Delta(S(a))=\zeta(S\ot S)\Delta(a)$ for all $a\in A$ (in an appropriate sense - see e.g.\ Proposition 3.7 and more comments in [VD-W3] for a correct formulation). Recall that we use $\zeta$ for the flip map. Moreover, the antipode satisfies the formulas 
$$\sum_{(a)} a_{(1)} S(a_{(2)}) a_{(3)}=a
\qquad\quad\text{and}\qquad\quad
\sum_{(a)} S(a_{(1)}) a_{(2)} S(a_{(3)})=S(a)
$$
for all $a$ in $A$. One has to multiply with an element of $A$, left or right, in order to be able to use the Sweedler notation, and so strictly speaking, the formulas hold in $M(A)$ (see also Remark 1.2 below). 
\snl
We  have the equalities
$$\align E(a\ot 1)&=\sum_{(a)} \Delta(a_{(1)})(1\ot S(a_{(2)})) \tag"(1.1)"\\
(1\ot a)E&=\sum_{(a)} (S(a_{(1)})\ot 1)\Delta(a_{(2)}) \tag"(1.2)"
\endalign$$
for all $a$. These equations are equivalent with
$$\align \Delta(c)(a\ot 1)&=\sum_{(a)} \Delta(ca_{(1)})(1\ot S(a_{(2)})) \tag"(1.3)"\\
(1\ot a)\Delta(b)&=\sum_{(a)} (S(a_{(1)})\ot 1)\Delta(a_{(2)}b) \tag"(1.4)"
\endalign$$
for all $a,b,c$. Observe that using the Sweedler notation in these formulas is just a matter of notation and nothing more. Indeed, the formula (1.3) above is a shorthand for the formula
$\Delta(c)(a\ot 1)=\sum_i\Delta(p_i)(1\ot S(q_i))$ where $\sum_i p_i\ot q_i=(c\ot 1)\Delta(a)$. This is true for all the formulas with the Sweedler notation we have here in this preliminary section. It illustrates a remark already made in the introduction.
\snl
In the regular case, we have that the antipode maps $A$ to itself and is bijective. In fact, this property of the antipode characterizes the {\it regular} weak multiplier Hopf algebras. 
\snl
In that case, we have the following counterparts of the formulas (1.1) and (1.2) above. We have
$$\align E(1 \ot a)&=\sum_{(a)} \Delta(a_{(2)})(S^{-1}(a_{(1)})\ot 1) \tag"(1.5)"\\
(a\ot 1)E&=\sum_{(a)} (1\ot S^{-1}(a_{(2)}))\Delta(a_{(1)}) \tag"(1.6)"
\endalign$$
for all $a$. Again these formulas can also be written as
$$\align \Delta(b)(1 \ot a)&=\sum_{(a)} \Delta(ba_{(2)})(S^{-1}(a_{(1)})\ot 1) \tag"(1.7)"\\
(a\ot 1)\Delta(c)&=\sum_{(a)} (1\ot S^{-1}(a_{(2)}))\Delta(a_{(1)}c) \tag"(1.8)"
\endalign$$
for all $a,b,c$.
\snl
Observe the following peculiarity in these formulas. The formulas (1.3) and (1.4) are true in the non-regular case but the expressions need not be in $A\ot A$. On the other hand, the formulas (1.7) and (1.8) only make sense in the regular case (as the inverse of $S$ is involved), while now the expressions are true in $A\ot A$.
\snl
We now make an important remark about the covering of the previous formulas.

\inspr{1.2} Remark \rm 
i) First rewrite the (images of the) canonical maps $T_1$ and $T_2$, and of $T_3$ and $T_4$ in the regular case, using the Sweedler notation, as
$$\align
&\Delta(a)(1\ot b)=\sum_{(a)}a_{(1)}\ot a_{(2)}b 
 \quad\text{}\quad
(c\ot 1)\Delta(a)=\sum_{(a)}ca_{(1)}\ot a_{(2)} \tag"(1.9)" \\
&(1\ot b)\Delta(a)=\sum_{(a)}a_{(1)}\ot ba_{(2)} 
\quad\text{}\quad 
\Delta(a)(c\ot 1)=\sum_{(a)}a_{(1)}c\ot a_{(2)} \tag"(1.10)" \\
\endalign$$
where $a,b,c\in A$. In all these four expressions, either $a_{(1)}$ is covered by $c$ and or $a_{(2)}$  by $b$. This is by the assumption put on the coproduct, requiring that the canonical maps have range in $A\ot A$.
\snl
ii)
Next consider the expressions
$$\align
&\sum_{(a)}a_{(1)}\ot S(a_{(2)})b \qquad \quad\text{and}\qquad\quad \sum_{(a)}cS(a_{(1)})\ot a_{(2)} \tag"(1.11)" \\
&\sum_{(a)}a_{(1)}\ot bS(a_{(2)}) \qquad \quad\text{and}\qquad\quad \sum_{(a)}S(a_{(1)})c\ot a_{(2)} \tag"(1.12)" \\
\endalign$$
where $a,b,c\in A$. In the first two formulas (1.11), we have a covering by the assumption that the generalized inverses $R_1$ and $R_2$ of the canonical maps exist as maps on $A\ot A$ with range in $A\ot A$ (see [VD-W3]). In the second pair of formulas (1.12), we have a good covering only in the regular case. It follows by considering the expressions in (1.9) and using that $S$ is a bijective anti-algebra map from $A$ to itself. In the regular case, we can also consider the above expressions with $S$ replaced by $S^{-1}$.
\snl
iii) If on the one hand, we first apply $S$ in the first or the second factor of the expressions in (1.9) and multiply and if on the other hand we simply apply multiplication on the expressions in (1.11), we get the four elements 
$$\align
&\sum_{(a)}S(a_{(1)})a_{(2)}b \qquad \quad\text{and}\qquad\quad \sum_{(a)}ca_{(1)}S(a_{(2)})  \\
&\sum_{(a)}a_{(1)}S(a_{(2)})b \qquad \quad\text{and}\qquad\quad \sum_{(a)}cS(a_{(1)})a_{(2)}  \\
\endalign$$
in $A$ for all $a,b,c\in A$. This is used to define the source and target maps in the next section (see Definition 2.1 in the next section).
\snl
iv) Now, we combine the coverings obtained in i) and ii). Consider e.g.\ the two expressions
$$\align
&\sum_{(a)} \Delta(a_{(1)})(1\ot S(a_{(2)})b) \tag"(1.13)"\\
&\sum_{(a)} (cS(a_{(1)})\ot 1)\Delta(a_{(2)}) \tag"(1.14)"
\endalign$$
where $a,b,c\in A$. The first expression (1.13) is obtained by applying the canonical map $T_1$ to the first of the two expressions in (1.11). So this gives an element in $A\ot A$ and we know that it is $E(a\ot b)$ as we can see from the formula (1.1). Similarly, the second expression (1.14) is obtained by applying the canonical map $T_2$ to the 
second of the two expressions in (1.11). We know that this is $(b\ot a)E$ as we see from the formula (1.2) above. Remark that $E(a\ot b)$ and $(b\ot a)E$ belong to $A\ot A$ because by assumption $E\in M(A\ot A)$, but that on the other hand, it is not obvious (as we see from the above arguments) that the expressions that we obtain for these elements belong to $A\ot A$.
\snl
v) Finally, as a consequence of the above statements, also the four expressions 
$$\align 
&\sum_{(a)} S(a_{(1)}) a_{(2)} S(a_{(3)})b \qquad\quad\text{and}\qquad\quad \sum_{(a)} ca_{(1)} S(a_{(2)}) a_{(3)}\\
&\sum_{(a)} cS(a_{(1)}) a_{(2)} S(a_{(3)}) \qquad\quad\text{and}\qquad\quad \sum_{(a)} a_{(1)} S(a_{(2)}) a_{(3)}b
\endalign$$
are well-defined in $A$ for all $a,b,c\in A$ (also in the non-regular case as $S:A\to M(A)$). This justifies a statement made earlier about the properties of the antipode. 
\hfill$\square$\einspr

And once again, in all these cases, the Sweedler notation is just used as a more transparent way to denote expressions. We refer to the coverings just to indicate how the formulas with the Sweedler notation can be rewritten without the use of it.
\snl
In the regular case, we also have many other nice formulas (see Section 4 in [VD-W3]. One of them is $(S\ot S)E=\zeta E$ (as expected because $E=\Delta(1)$). Other formulas that we will use, will be recalled later. In any case, they are all found in [VD-W3] and we refer to this paper for details.
\nl
\nl

\bf 2. The source and target algebras \rm
\nl
As in the previous section we consider a weak multiplier Hopf algebra $(A,\Delta)$. In general, we do not assume that it is regular. In the regular case, nicer results can be obtained, but we try to push the theory as far as possible in the general case.
\snl
We first recall the definition of the {\it source} and {\it target maps} $\varepsilon_s:A\to M(A)$ and $\varepsilon_t:A\to M(A)$ and prove the first properties. We show among other things that the images are non-degenerate subalgebras of $M(A)$, sitting nicely in $M(A)$ so that also their multiplier algebras can be considered as subalgebras of $M(A)$. 
\snl
The source and target maps, together with their images, have been considered already in [VD-W3] and a few properties were proven, mainly for the purpose of studying the antipode. In this paper, we will continue this study. 
\snl
Remark that in this paper, as we mentioned already in the introduction, {\it we will define the source and target algebras as the images of the source and target maps} (see Notation 2.9). We will explain later why we do this.
\snl
We will also study the behavior of the antipode $S$ on the source and target algebras. Recall that $S$ is an anti-homomorphism from $A$ to $M(A)$. It is non-degenerate in the sense that $S(A)A=A$ and $AS(A)=A$ (see Proposition 3.6 in [VD-W3]). Therefore, as a consequence of a general property mentioned already in the introduction (see also the appendix of [VD1]), it has a unique extension to a unital anti-homomorphism from $M(A)$ to itself. 
\snl
We consider the canonical idempotent $E$ in $M(A\ot A)$ as reviewed in the previous section and we use that the coproduct $\Delta$ can be extended to the multiplier algebra as we have mentioned earlier.
We will show that $E$ is a separability idempotent as studied in [VD4.v2].
\nl
\it The source and target algebras $B$ and $C$ \rm
\nl
We first consider the {\it source} and {\it target maps} $\varepsilon_s:A\to M(A)$ and $\varepsilon_t:A\to M(A)$. Recall Definition 3.1 from [VD-W3].

\inspr{2.1} Definition \rm
For $a\in A$ we define 
$$\varepsilon_s(a)=\sum_{(a)}S(a_{(1)})a_{(2)}
\qquad\quad\text{and}\qquad\quad
\varepsilon_t(a)=\sum_{(a)}a_{(1)}S(a_{(2)})$$
where $S$ is the antipode. The map $\varepsilon_s$ is called the source map and the map $\varepsilon_t$ is the target map. 
\hfill $\square$ \einspr

We have seen in Remark 1.2.iii that these maps have well-defined values in the multiplier algebra $M(A)$.
\snl
We will show that the images of the source and target maps are subalgebras of $M(A)$. Before we can do this, we need some elementary properties, also important for the further study of these subalgebras.
\snl
First we have that the range of $\varepsilon_s$ coincides with the left leg of $E$ and that the range of $\varepsilon_t$ is the right leg of $E$. These statements are made precise in the following proposition.

\inspr{2.2} Proposition \rm
The range $\varepsilon_s(A)$ of the source map is spanned by elements of the form \newline
$(\iota\ot\omega(a\,\cdot\,b))E$
where $a,b\in A$ and where $\omega$ is a linear functional on $A$. Similarly the range $\varepsilon_t(A)$ of the target map is spanned by elements of the form
$(\omega(c\,\cdot\,a)\ot\iota)E$
where $a,c\in A$ and with $\omega$ a linear functional on $A$.

\snl\bf Proof\rm:
By formula (1.2) in Section 1 we get for $a,b\in A$ that
$$(1\ot a)E(1 \ot b)=\sum_{(a)} S(a_{(1)})a_{(2)} \ot a_{(3)}b$$
and this belongs to $\varepsilon_s(A)\ot A$. We can apply a linear functional $\omega$ on the second leg and we see that $(\iota\ot\omega(a\,\cdot\,b))E$ is well-defined and belongs to $\varepsilon_s(A)$. The fullness of $\Delta$ guarantees that any element of $A$ is a sum of elements of the form $$(\iota\ot\omega)(\Delta(a)(1\ot b))$$
where $a,b\in A$ and where $\omega$ is a linear functional (see Section 1). Hence it follows that $\varepsilon_s(A)$ is spanned by elements as in the formulation of the proposition. Similarly for the range $\varepsilon_t(A)$ of the target map.
\hfill $\square$\einspr

Because $E\ot 1$ and $1\ot E$ commute, it follows that $\varepsilon_s(a)$ and $\varepsilon_t(b)$ will commute in $M(A)$ for all $a,b\in A$.
\snl
Also the following is an easy consequence of the previous result. The formulas in the proposition make sense in the multiplier algebra $M(A\ot A)$.

\inspr{2.3} Proposition \rm
$$\Delta(x)=(x\ot 1)E=E(x\ot 1)
\quad\qquad\text{and}\quad\qquad
\Delta(y)=E(1\ot y)=(1\ot y)E$$
for $x\in \varepsilon_t(A)$ and $y\in \varepsilon_s(A)$. 
 
\bf \snl Proof\rm:
Simply apply the appropriate linear functionals on the first, respectively the third factor of the equations
$$\align
(\iota\ot \Delta)E&=(E\ot 1)(1\ot E)=(1\ot E)(E\ot 1)\\
(\Delta\ot\iota)E&=(E\ot 1)(1\ot E)=(1\ot E)(E\ot 1).
\endalign$$
\vskip -0.5 cm
\hfill $\square$ \einspr

The result above is the motivation for the following lemma.

\inspr{2.4} Lemma \rm
For an element $x\in M(A)$, the following are equivalent:
\snl
\hbox{\hskip 0.3 cm} i) $\Delta(x)=(x\ot 1)E$, \newline
\hbox{\hskip 0.3 cm} ii) $\Delta(x)=E(x\ot 1)$. 
\snl
Similarly, for an element $y\in M(A)$, the following are equivalent:
\snl
\hbox{\hskip 0.3 cm} i) $\Delta(y)=E(1\ot y)$, \newline
\hbox{\hskip 0.3 cm} ii) $\Delta(y)=(1\ot y)E$.

\bf \snl Proof\rm: First let $x\in M(A)$ and assume that $\Delta(x)=(x\ot 1)E$. Take any 
 $y\in \varepsilon_s(A)$. Then
$$\Delta(xy)=\Delta(x)\Delta(y)= (x\ot 1)E\Delta(y)=(x\ot 1)\Delta(y)=(x\ot y)E.$$
We have used that $\Delta(y)=(1\ot y)E$, proven in the previous proposition for elements $y$ in $\varepsilon_s(A)$.
On the other hand
$$\Delta(yx)=\Delta(y)\Delta(x)=(1\ot y)E\Delta(x)=(1\ot y)\Delta(x)=(x\ot y)E$$
and we see that $\Delta(xy)=\Delta(yx)$. Multiply with $\Delta(a)$ for any $a\in A$ and apply the counit. This will give $xya=yxa$ and because this is true for all $a$, we have $xy=yx$.
\snl
Because this result is true for all elements $y$ in the left leg of $E$, as a consequence we find that $(x\ot 1)E=E(x\ot 1)$ and hence also $\Delta(x)=E(x\ot 1)$.
\snl 
Similarly if $\Delta(x)=E(x\ot 1)$ also $\Delta(x)=(x\ot 1)E$ will be true. This proves the equivalence of i) and ii) in the first part of the lemma.
\snl
The second part is proven in a completely similar way.
\hfill $\square$ \einspr 

We arrive at the following notation.

\inspr{2.5} Notation \rm
We will denote by $A_s$ the set of elements $y\in M(A)$ satisfying $\Delta(y)=E(1\ot y)$ and by $A_t$ the set of elements $x\in M(A)$ satisfying $\Delta(x)=(x\ot 1)E$.
\hfill $\square$ \einspr 

The following is an immediate consequence of the lemma.

\inspr{2.6} Proposition \rm The sets $A_s$ and $A_t$ are commuting subalgebras of $M(A)$. 

\snl\bf Proof\rm: 
It is immediately clear from the definitions that these sets are subalgebras of $M(A)$. Moreover, if $x\in A_t$ and $y\in A_s$, we have as in the first part of the proof of the lemma
$$
\Delta(xy)=\Delta(x)\Delta(y)=(x\ot y)E
\tussen
\Delta(yx)=\Delta(y)\Delta(x)=(x\ot y)E
$$
where we have used the two equivalences of i) and ii) in the lemma. Hence $\Delta(xy)=\Delta(yx)$ and as before, $xy=yx$.
\hfill $\square$ \einspr 

From Proposition 2.3, we know that 
$\varepsilon_s(A)\subseteq A_s$ and $\varepsilon_t(A)\subseteq A_t$.
However, we can now prove more.

\inspr{2.7} Proposition \rm
Assume that $x\in A_t$. Then for all $a\in A$ we have
$\varepsilon_t(xa)=x\varepsilon_t(a)$ and $\varepsilon_s(ax)=S(x)\varepsilon_s(a)$. Similarly, if $y\in A_s$ we have $\varepsilon_s(ay)=\varepsilon_s(a)y$
and $\varepsilon_t(ya)=\varepsilon_t(a)S(y)$ for all $a\in A$.

\snl \bf Proof\rm:
Take $x\in M(A)$ and assume that $\Delta(x)=(x\ot 1)E$. Let $a\in A$. Then
$\Delta(xa)=(x\ot 1)\Delta(a)$ and if we apply $m(\iota \ot S)$ where $m$ is multiplication, we find 
$\varepsilon_t(xa)=x\varepsilon_t(a)$. By Lemma 2.4, we know that also $\Delta(ax)=\Delta(a)(x\ot 1)$ and now we apply $m(S\ot \iota)$ to find $\varepsilon_s(ax)=S(x)\varepsilon_s(a)$. This proves the first part of the proposition.
\snl
The second part is proven in a completely similar way.
\hfill$\square$\einspr

Using techniques as above, we find other formulas of this type but we will not need these.
\snl
The result above has a few obvious, but important consequences.

\inspr{2.8} Proposition \rm 
i) The sets $\varepsilon_s(A)$ and $\varepsilon_t(A)$ are subalgebras.\newline
ii) The algebra $\varepsilon_s(A)$ is a right ideal of $A_s$ and $\varepsilon_t(A)$ is a left ideal of $A_t$.
\hfill $\square$ \einspr

Remark that the algebras $A_s$ and $A_t$ contain the identity of $M(A)$. This is not the case in general for the subalgebras $\varepsilon_s(A)$ and $\varepsilon_t(A)$. It is also not clear if, again in general, $\varepsilon_s(A)$ is also a left ideal of $A_s$ and if $\varepsilon_t(A)$ is also a right ideal of $A_t$. All of this is related with the behavior of the antipode on these algebras (as we can see already from formulas in Proposition 2.7). In a subsequent item, we investigate this further.
\snl
First we look at the multiplier algebras of the images of the source and the target maps.
\nl
\it The multiplier algebras of the source and target algebras \rm
\nl
We introduce the following notation and terminology. As mentioned already in the introduction, the terminology is different from the one originally used in [VD-W3], see further.

\inspr{2.9} Notation \rm
In what follows, we will denote the algebra $\varepsilon_s(A)$ by $B$ and $\varepsilon_t(A)$ by $C$. We will call $B$ the {\it source algebra} and $C$ the {\it target algebra}.
\hfill $\square$ \einspr

Recall that we do not expect these algebras to be unital. And we are interested in their multiplier algebras, if they exist.
\snl
We begin with some module properties giving more information about these algebras $B$ and $C$ and how they sit in $M(A)$. 

\iinspr{2.10} Proposition \rm
We have 
$$\align A&=AB\qquad\qquad\text{and}\qquad\qquad A=CA\\
A&=BA\qquad\qquad\text{and}\qquad\qquad A=AC.
\endalign$$

\bf Proof\rm:
We know that 
$$ba=\sum_{(a)} ba_{(1)} S(a_{(2)}) a_{(3)}=\sum_{(a)} ba_{(1)} \varepsilon_s(a_{(2)})$$
for all $a,b$. The right hand side is in $A\varepsilon_s(A)$ and because $A^2=A$ we find that $A=A\varepsilon_s(A)$. Similarly, from the formula
$$ab=\sum_{(a)} a_{(1)} S(a_{(2)}) a_{(3)}b=\sum_{(a)} \varepsilon_t(a_{(1)}) a_{(2)}b$$
for all $a,b$, we get $A=\varepsilon_t(A)A$.
\snl
If on the other hand, we start with the formula 
$$bS(a)=\sum_{(a)} bS(a_{(1)}) a_{(2)} S(a_{(3)})=\sum_{(a)} bS(a_{(1)}) \varepsilon_t(a_{(2)})$$
for all $a,b$, we find that $AS(A)$ is contained in $A\varepsilon_t(A)$ (recall Remark 1.2.ii in Section 1). Now, in Proposition 3.6 of [VD-W3], we have shown that $AS(A)=A$ and so we get also $A=A\varepsilon_t(A)$. Similarly $A=\varepsilon_s(A)A$.
\hfill$\square$\einspr  

The results above say that $A$ as a $B$-bimodule and as a $C$-bimodule is unital. If we combine the above result with the property in Proposition 2.7 we get the following.

\iinspr{2.11} Corollary \rm 
The algebras $B$ and $C$ are idempotent. 
\hfill $\square$ \einspr

Indeed, for all $a,b$ we have e.g.\ $\varepsilon_s(a\varepsilon_s(b))=\varepsilon_s(a)\varepsilon_s(b)$. Similarly for $\varepsilon_t(A)$.
\snl
Later, we will see that the algebras $B$ and $C$ have  local units. This implies that the bimodules are also non-degenerate. In fact, this  already follows by a more general argument, which is part of the following, also more general result.

\iinspr{2.12} Lemma \rm
Let $R$ be a subalgebra of $M(A)$. Multiplication makes $A$ into a $R$-bimodule. Assume that this module is unital. Then it is also a non-degenerate bimodule. The algebra $R$ is a non-degenerate algebra and the embedding of $R$ in $M(A)$ extends uniquely to an embedding of the multiplier algebra $M(R)$ of $R$ in $M(A)$. Moreover we have, considering $M(R)$ as sitting inside $M(A)$,
$$M(R)=\{x\in M(A) \mid xr\in R \text{ and } rx\in R \text{ for all } r\in R \}.\tag"(2.1)"$$

\snl\bf Proof\rm:
We first show that the module is non-degenerate. Take any $a\in A$ and assume that $ra=0$ for all $r\in R$. Then $a'ra=0$ for all $a'\in A$ and $r\in R$. Because we assume that $AR=A$, it follows that also $a'a=0$ for all $a'\in A$. Then $a=0$. Similarly on the other side. We  get in that $A$ is a non-degenerate $R$-bimodule. 
\snl
We also claim that  $R$ is a non-degenerate subalgebra of $M(A)$. To show this assume that $r\in R$ and that $rs=0$ for all $s\in R$. Multiply with an element $a\in A$ from the right and use that $RA=A$. This implies that $ra=0$ for all $a\in A$. Then $r=0$. Similarly on the other side. So the algebra $R$ is non-degenerate and we can consider its multiplier algebra $M(R)$.
\snl
As $A$ is assumed to be a unital $R$-bimodule, we have $RA=A$ and $AR=A$. So the embedding $j:R \to M(A)$ is a non-degenerate homomorphism and a standard result implies that it extends uniquely to a unital homomorphism $\overline j: M(R) \to M(A)$. It is not hard to show that in this case, this extension is still an embedding. Because obviously for any $x\in M(R)$ we have $xr\in R$ and $rx\in R$ for all $r\in R$, we find one inclusion of the statement (2.1). The other inclusion is proven by using again that the $R$-bimodule $A$ is unital.
\hfill$\square$\einspr

We can apply this lemma and we obtain the following. Recall that we use $B$ to denote the algebra $\varepsilon_s(A)$ and $C$ for $\varepsilon_t(A)$ (cf.\ Notation 2.9).

\iinspr{2.13} Proposition \rm 
The algebras  $B$ and $C$ are non-degenerate and idempotent. Their multiplier algebras $M(B)$ and $M(C)$ embed in $M(A)$. An element $x\in M(C)$ still satisfies
$$\Delta(x)=(x\ot 1)E=E(x\ot 1)$$
while 
$$\Delta(y)=E(1\ot y)=(1\ot y)E$$
is still true for elements $y$ of $M(B)$. So $M(B)\subseteq A_s$ and $M(C)\subseteq A_t$.

\snl\bf Proof\rm:
The conditions in Lemma 2.12 are fulfilled for the subalgebras $B$ and $C$ as shown in the Propositions 2.10. Therefore, $B$ and $C$ are non-degenerate algebras and they sit in $M(A)$ in such a way that the embeddings $B\subseteq M(A)$ and $C\subseteq M(A)$ extend to embedding of their multiplier algebras $M(B)$ and $M(C)$.
\snl
We have already explained that the algebras $B$ and $C$ are idempotent.
There are various ways to prove that we still have the embeddings $M(B)\subseteq A_s$ and $M(C)\subseteq A_t$. Take e.g. $m\in M(C)$, $x\in C$ and $a\in A$. Then
$$\Delta(mxa)=(mx\ot 1)\Delta(a)=(m\ot 1)\Delta(xa).$$
As $CA=A$, it follows that $\Delta(ma)=(m\ot 1)\Delta(a)$ for all $a\in A$ and hence $\Delta(m)=(m\ot 1)E$. Similar arguments are used for the other equations.
\hfill$\square$\einspr

In the next item of this section, we study the behavior of the antipode on the algebras $B$ and $C$.

\snl
\it The antipode on the source and target algebras \rm
\nl
We begin with the following result about the subalgebras $A_s$ and $A_t$ of $M(A)$. Recall that we can extend the antipode $S$ to a unital anti-homomorphism from $M(A)$ to itself. 

\iinspr{2.14} Proposition \rm 
i) If $x,y\in M(A)$ and $(1\ot x)E=(y\ot 1)E$, then $x\in A_t$ and $y\in A_s$. \newline
ii) If $x,y\in M(A)$ and $E(1\ot x)=E(y\ot 1)$, then $x\in A_t$ and $y\in A_s$. \newline
iii) If $x\in A_t$ then $S(x)\in A_s$ and $(1\ot x)E=(S(x)\ot 1)E$.\newline
iv) If $y\in A_s$ then $S(y)\in A_t$ and $E(y\ot 1)=E(1\ot S(y))$.

\snl\bf Proof\rm:
i) Assume $x,y\in M(A)$ and that $(1\ot x)E=(y\ot 1)E$.  If we apply $\iota\ot\Delta$ to this equation, we find
$$\align 
(1\ot \Delta(x))(E\ot 1)
&=(y\ot 1\ot 1)(E\ot 1)(1\ot E)\\
&=(1\ot x\ot 1)(E\ot 1)(1\ot E)\\
&=(1\ot x\ot 1)(1\ot E)(E\ot 1).
\endalign$$
Now we use the property that $(1\ot a)E=0$ implies that $a=0$ (see Lemma 1.1 in Section 1). This will eventually give $\Delta(x)=(x\ot 1)E$. This proves that $x\in A_t$. If we apply $\Delta\ot \iota$ instead, we  obtain that $y\in A_s$. 
\snl
ii) The second property is proven in completely the same way.
\snl
iii) Let $x\in A_t$ so that $\Delta(x)=E(x\ot 1)$. Then for all $a\in A$ we have $\Delta(ax)=\Delta(a)(x\ot 1)$ and so
$$\align (1\ot ax)E
&=\sum_{(ax)} (S((ax)_{(1)}) \ot 1)\Delta((ax)_{(2)})\\
&=\sum_{(a)} (S(a_{(1)}x) \ot 1)\Delta(a_{(2)})\\
&=\sum_{(a)} (S(x)S(a_{(1)}) \ot 1)\Delta(a_{(2)})\\
&=(S(x)\ot a)E
\endalign$$
This implies $(1\ot x)E=(S(x)\ot 1)E$. It follows from i) that $S(x)\in A_s$.
\snl
iv) Similarly we get $S(y)\in A_t$ when $y\in A_s$ and $E(y\ot 1)=E(1\ot S(y))$.
\hfill $\square$ \einspr

Remark that it follows that $S$ is injective on $A_s$ and on $A_t$. However, it does not imply that these maps are surjective in the general case.
\nl
We now investigate the maps $S_B:B\to M(A)$ and $S_C:C\to M(A)$ that we obtain by restricting (the extension of) the antipode to the subalgebras $B$ and $C$ of $M(A)$. As a special case of the equations above, we have
$$(1\ot x)E=(S_C(x)\ot 1)E
\qquad\quad\text{and}\qquad\quad
E(y\ot 1)=E(1\ot S_B(y))$$
for $x\in B$ and $y\in C$. 
In particular, we know already that $S_B:B\to A_t$ and $S_C:C\to A_s$. In the next proposition, we get a stronger result.

\iinspr{2.15} Proposition \rm
The map $S_B$ is a non-degenerate anti-homomorphism from $B$ to $M(C)$ and the map $S_C$ is a non-degenerate anti-homomorphism from $C$ to $M(B)$. Both maps are injective.

\snl\bf Proof\rm:
i) Take $x\in C$. Then $x\in A_t$ and from Proposition 2.7 we know that $\varepsilon_s(ax)=S(x)\varepsilon_s(a)$ for all $a$. Because now also $\varepsilon_s(aS(x))=\varepsilon_s(a)S(x)$ for all $a$, we see that $S(x)\in M(B)$. Similarly $S(y)\in M(C)$ when $y\in B$. It follows that $S_C$ is an anti-homomorphism from $C$ to $M(B)$ and that $S_B$ is an anti-homomorphism of $B$ to $M(C)$.
\snl
ii) As $BA=A$ and $\varepsilon_t(ya)=\varepsilon_t(a)S(y)$ for $y\in B$, we see that $CS(B)=C$. On the other hand we have 
$$A=S(A)A=S(AB)A=S(B)S(A)A=S(B)A$$
and because $\varepsilon_t(S(y)a)=S(y)\varepsilon_t(a)$ for $y\in B$ we see that also $S(B)C=C$. 
\snl
Hence $S_B:B\to M(C)$ and $S_C:C\to M(B)$ are non-degenerate anti-homomor\-phisms.
\hfill $\square$ \einspr

From the general theory, we know that $S_B$ and $S_C$ have unique extensions to unital anti-homomor\-phism from $M(B)$ to $M(C)$ and from $M(C)$ to $M(B)$ respectively. These extensions are still the restrictions of the antipode $S$ to the multiplier algebras $M(B)$ and $M(C)$ respectively. 

\snl
In the {\it regular case}, we have the following stronger results.

\iinspr{2.16} Proposition \rm
In the case of a regular weak multiplier Hopf algebra, we have that $S_B$ is an anti-isomorphism from $B$ tot $C$ and $S_C$ is an anti-isomorphism from $C$ to $B$.
The multiplier algebras $M(B)$ and $M(C)$ are respectively equal to the algebras $A_s$ and $A_t$ as defined in Notation 2.5.

\snl\bf Proof\rm:
We can use e.g.\ that $(S\ot S)E=\zeta E$ in the case of a regular weak multiplier Hopf algebra (see Proposition 4.4 in [VD-W3]). As $B$ is the left leg of $E$ and $C$ is the right leg of $E$, we find that $S$ maps $B$ to $C$ and $C$ to $B$. It also follows that these maps are surjective. As we know already that they are also injective, we find the first statement of the proposition.
\snl
The equation $(S\ot S)E=\zeta E$ also implies that $S$ maps $A_s$ to $A_t$ and vice versa. In Proposition 2.7 we have shown that 
$$\varepsilon_s(ay)=\varepsilon_s(a)y
\tussenen
\varepsilon_s(ax)=S(x)\varepsilon_s(a)$$
for $y\in A_s$ and $x\in A_t$. It follows that the algebra $B$, the image $\varepsilon_s(A)$, is a two-sided ideal of $A_s$. And because we know already that the $M(B)\subseteq A_s$, it follows that $M(B)=A_s$.
\snl
Similarly we have $M(C)=A_t$.
\hfill $\square$\einspr

It is not completely clear what the situation is in the non-regular case. We have Proposition 2.15 saying that $S_B$ embeds $B$ in $M(C)$ and Proposition 2.13 saying that $M(C)$ is a subalgebra of $A_t$. Similarly $S_C$ embeds $C$ in $M(B)$ and $M(B)$ is a subalgebra of $A_s$. 
\snl
For this reason, we have changed our terminology and are now calling the algebras $B$ and $C$, the images of the source and target maps respectively, the source and target algebras. In an earlier version of this paper [VD-W4.v1] we have used these terms for $A_s$ and $A_t$ instead. This was motivated by the fact that, in the regular case, the can be identified with the multiplier algebras of $B$ and $C$ respectively. But this is not sure in the non-regular case that we are investigating in greater detail in this version of the paper. 
\nl
\it The canonical idempotent $E$ as a separability idempotent in $M(B\ot C)$ \rm
\nl
We have the algebras $B$ and $C$. They are non-degenerate and idempotent. The algebra $B$ is the left leg of $E$ and the algebra $C$ is the right leg of $E$, in an appropriate sense, see Proposition 2.2. And because $E$ is a multiplier of $A\ot A$, we can expect that it is also a multiplier of $B\ot C$. This turns out to be the case. Moreover it is a separability idempotent as defined and studied in [VD4.v2]. This is what we show next.
\snl
The first step is the following result.

\iinspr{2.17} Lemma \rm We have 
$$E(1\ot a)\in B\ot A
\tussenen
(a\ot 1)E\in A\ot C$$
for all $a\in A$.

\snl\bf Proof\rm:
For all $a$ in $A$ we can define a left multiplier $\varepsilon'_s(a)$ of $A$ by the formula
$$\varepsilon'_s(a)b=(\iota\ot\varepsilon)(E(b\ot a))$$
where $b$ is in $A$. We will see later why we use this notation. 
\snl
Fix two elements $a,a'$ in $A$. Write 
$$\sum_{(a)} \varepsilon'_s(a_{(1)})\ot a_{(2)}a'=\sum_i t_i\ot q_i$$
where $t_i$ is a left multiplier of $A$ and $q_i\in A$. Assume that the $(q_i)$ are linearly independent. 
\snl 
For all $b$ in $A$ we find
$$\align
\sum_{(a)} \varepsilon'_s(a_{(1)})b\ot a_{(2)}a'
&=\sum_{(a)} (\iota\ot\varepsilon\ot\iota)((E\ot 1)(b\ot a_{(1)}\ot a_{(2)}a'))\\
&=(\iota\ot\varepsilon\ot\iota)((E\ot 1)(1\ot E)(b\ot \Delta(a)(1\ot a')))\\
&=(\iota\ot\varepsilon\ot\iota)((\iota\ot\Delta)(E(b\ot a))(1\ot 1\ot a'))\\
&=E(b\ot a)(1\ot a')=E(b\ot aa').\endalign
$$
Therefore $E(b\ot aa')=\sum_i t_ib\ot q_i$
for all $b\in A$.
\snl
On the other hand, for all $c\in A$, we have also
$$\align 
(1\ot c)E(1\ot aa')
&=\sum_{(c)}S(c_{(1)})c_{(2)}\ot c_{(3)}aa'\\
&=\sum_{(c)}\varepsilon_s(c_{(1)})\ot c_{(2)}aa'
\endalign$$
and this belongs to $B\ot A$.
\snl
If we combine this with the previous formulas, we find $\sum_i t_i\ot cq_i\in B\ot A$ for all $c\in A$. Now let $\omega$ be a linear functional on the space $L(B)$ of left multipliers of $A$ that vanishes on elements in $B$. We find $\sum_i \omega(t_i)cq_i=0$ for all $c$ in $A$. By non-degeneracy of the product in $A$ and because the elements $(q_i)$ are linearly independent, it follows that $\omega(t_i)=0$ for all $i$. Hence $t_i$ is in $B$ for all $i$ and we find that $E(1\ot aa')\in B\ot A$. Because $A$ is idempotent, we get $E(1\ot A)\subseteq B\ot A$.
\snl
In a completely similar way we can prove that also $(A\ot 1)E\in A\ot C$. This proves the lemma.
\hfill $\square$ \einspr

From the proof we see that $\sum_{(a)} \varepsilon'_s(a_{(1)})\ot a_{(2)}a'\in B\ot A$ and from the fullness of the coproduct, it follows that $\varepsilon'(a)\in B$ for all $a\in A$. This in turn follows of course also from the property that $E(1\ot A)\subseteq B\ot A$. 
\snl
We will give more comments on this result later. First we use the lemma to prove the following main result.

\iinspr{2.18} Theorem \rm 
The canonical idempotent of a weak multiplier Hopf algebra is a separability idempotent in $M(B\ot C)$ where $B$ and $C$ are the source and target algebras.

\snl \bf Proof\rm:
i) By the lemma, we find that $E(1\ot a)$ belongs to $B\ot A$. We therefore can apply $\varepsilon_t$ on the second leg of this expression. We know that the second leg of $E$ belongs to $\varepsilon_t(A)$ and this is a subalgebra of $A_t$. In Proposition 2.7 we have shown that $\varepsilon_t(xa)=x\varepsilon_t(a)$ for all $x\in A_t$. Therefore $(\iota\ot\varepsilon_t)(E(1\ot a))=E(1\ot \varepsilon_t(a))$. We conclude that $E(1\ot \varepsilon_t(a))\in B\ot C$ for all $a$ and so $E(1\ot C)\subseteq B\ot C$. 
\snl
In a completely similar way, we find that $(B\ot 1)E\subseteq B\ot C$. It follows not only that $E\in M(B\ot C)$, but also that it satisfies the first requirements for a separability idempotent (see Section 1 of [VD4.v2]).
\snl
ii) We will now show that $E$ is full in the sense of Definition 1.1 of [VD4.v2]. For this, assume that $V$ is a subspace of $B$ so that $E(1\ot x)\subseteq V\ot C$ for all $x\in C$. Then $(1\ot b)E(1\ot xa)\in V\ot A$ for all $a,b\in A$ and $x\in C$. In Proposition 2.10 we showed that $CA=A$ and in Proposition 2.2 that $B$  is spanned by elements of the form
$(\iota\ot\omega(a\,\cdot\,b))E$
where $a,b\in A$ and where $\omega$ is a linear functional on $A$. Then we must have $V=B$ proving that the left leg of $E$ (as an idempotent in $M(B\ot C)$) is still all of $B$. Similarly for the right leg. Hence $E$ is full.
\snl
iii) Finally, we know already from Proposition 2.15 that the antipode is a non-degenerate anti-homomorphism from $B$ to $M(C)$ as well as a non-degenerate anti-homomorphism from $C$ to $M(B)$. As in Proposition 2.14 they satisfy 
$$(1\ot x)E=(S(x)\ot 1)E
\tussenen 
E(y\ot 1)=E(1\ot S(y))$$
when $x\in C$ and $y\in B$. This is the final requirement in Definition 1.4 of [VD4.v2] and shows that $E$ is a separability idempotent in $M(B\ot C)$. This completes the proof.
\hfill $\square$ \einspr

Remark that in item iii) of the proof above, we find $E(y\ot 1)=E(1\ot S(y))$ for all $y\in B$. Then $E(1\ot S(y)x)=E(y\ot x)$ for all $x\in C$ and $y\in B$. From the fact that $E\in M(B\ot C)$ and that $S$ is a non-degenerate anti-homomorphism from $B$ to $M(C)$, it would also follow that $E(1\ot C)\subseteq B\ot C$.
\snl
In the {\it regular case}, we have the following expected result.

\iinspr{2.19} Proposition \rm
If the weak multiplier Hopf algebra is regular, then $E$ is a regular separability idempotent.

\snl \bf Proof\rm:
There are different ways to prove this. If we start with the definition of regularity for a weak multiplier Hopf algebra (as e.g.\ in Definition 4.1 of [VD-W3]), then we assume that $(A,\Delta^{\text{cop}})$ also satisfies the axioms of a weak multiplier Hopf algebra. The canonical idempotent now is $\zeta E$ where $E$ is the canonical idempotent of the original weak multiplier Hopf algebra. Remember that $\zeta$ is the flip map on $A\ot A$ and extended to $M(A\ot A)$. 
\snl
Because $B$ and $C$ are the left and the right legs of $E$, we get that $C$ and $B$ are the left and the right legs of $\zeta E$. Applying Theorem 2.18 to the new weak multiplier Hopf algebra $(A,\Delta^{\text{cop}})$, we obtain that $\zeta E$ is a separability idempotent in $M(C\ot B)$. Then $E$ is indeed a regular separability idempotent by the very definition of regularity for a separability idempotent (see Definition 2.4 of [VD4.v2]).
\hfill $\square$ \einspr

In an earlier version of this paper (reference [VD-W4.v1]), we only considered regular weak multiplier Hopf algebras and this result was obtained already, see Section 2 in [VD-W4.v1]).
\nl
Let us now consider some of the results we have proven for general and regular separability idempotents in [VD4.v2] and see what they give in the case of the canonical idempotent of a weak multiplier Hopf algebra. Recall the distinguished linear functionals $\varphi_B$ and $\varphi_C$ on $B$ and $C$ respectively, defined and characterized by the formulas
$$(\varphi_B\ot\iota)E=1
\tussenen
(\iota\ot\varphi_C)E=1;$$
see Proposition 1.9 in [VD4.v2].

\iinspr{2.20} Proposition \rm
The distinguished linear functionals $\varphi_B$ and $\varphi_C$, obtained for the separability idempotent $E$, satisfy
$$\varphi_B(\varepsilon_s(a))=\varepsilon(a)
\tussenen
\varphi_C(\varepsilon_t(a))=\varepsilon(a)
$$
for all $a\in A$.

\snl \bf Proof\rm: 
We have the formula
$$(1\ot a)E(1\ot b)=\sum_{(a)}\varepsilon_s(a_{(1)})\ot a_{(2)}b$$
for all $a,b\in A$ (see e.g.\ in the proof of Lemma 2.17). If we apply $\varphi_B$ on the first factor, we obtain
$$\varphi_B(\varepsilon_s(a_{(1)}))a_{(2)}b=ab.$$
If we apply a linear functional $\omega$ we find $\varphi_B(\varepsilon_s(a'))=\omega(ab)$ with  
$$a'=(\iota\ot\omega)(\Delta(a)(1\ot b).$$ 
Because $\varepsilon(a')=\omega(ab)$ we see that $\varphi_B(\varepsilon_s(a'))=\varepsilon(a')$. By the fullness of the coproduct, any element of $A$ is of the form $(\iota\ot\omega)(\Delta(a)(1\ot b)$. This proves the first formula of this proposition. The other one is proven in a similar way.
\hfill $\square$ \einspr 
\snl
\it Existence of local units \rm
\nl
From the general theory of (possibly non-regular) separability idempotents, we know that there exist local units (cf. Proposition 1.10 in [VD4.v2]). As a consequence we get the following result. 

\iinspr{2.21} Proposition \rm
The algebra $A$ has local units.

\snl\bf Proof\rm:
Let $a\in A$ and assume that $\omega$ is a linear functional on $A$ so that $\omega(ba)=0$ for all $b\in A$. Then 
$$(\iota\ot\omega)((1\ot b)(\iota\ot S)((c\ot 1)\Delta(p))(1\ot a))=0$$
for all $b,c,p\in A$. We use that $(c\ot 1)\Delta(p)\in A\ot A$. We know that $((\iota\ot S)\Delta(p))(1\ot a)$ belongs to $A\ot A$. Therefore, we can cancel $c$ in the above equation and get
$$(\iota\ot\omega)((1\ot b)((\iota\ot S)\Delta(p))(1\ot a))=0.$$
Write $((\iota\ot S)\Delta(p))(1\ot a)$ as $\sum_i p_i\ot q_i$ and assume that the elements $(p_i)$ are linearly independent. We find $\omega(bq_i)=0$ for all $i$ and all $b\in A$. Replace $b$ by $p_i$ and take the sum over $i$. Because $\sum p_iq_i=\varepsilon_t(p)a$ we get $\omega(\varepsilon_t(p)a)=0$ for all $p\in A$. This means that $\omega(xa)=0$ for all $x\in C$.
\snl
We know that $A=CA$ and because we have left local units in $C$, there exists an element $x\in C$ so that $xa=a$. Then we see that $\omega(a)=0$. This means that $a\in Aa$ and we know that this implies that $A$ has left local units. In a similar way, we find that $A$ also has right local units. This completes the proof. 
\hfill $\square$ \einspr

We see in the proof that we only need that $B$ has right local units and that $C$ has left local units. These results have a more easy proof in [VD4.v2].
\snl
Recall also that in earlier work on weak multiplier Hopf algebras, the existence of local units was only obtained in the case of a regular weak multiplier Hopf algebra, see Proposition 4.9 in [VD-W3].
\nl
We finish this section with a couple of remarks.

\iinspr{2.22} Remarks \rm
i) As we see from the proof of Lemma 2.17 and from earlier arguments, we find that 
$(\iota\ot\varepsilon)((1\ot a)E)=\varepsilon_s(a)$ when $a\in A$. The formula makes sense as an equality of left multipliers of $A$. Remark that we do not expect $(1\ot a)E$ to belong to $B\ot A$. Similarly, we find 
$(\varepsilon\ot\iota)(E(a\ot 1))=\varepsilon_t(a)$ for $a$ in $A$, now as right multipliers of $A$. Again we do not expect $E(A\ot 1)\subseteq A\ot C$.
\snl
ii) On the other hand, we do have $E(1\ot A)\subseteq B\ot A$ and $(A\ot 1)E\subseteq A\ot C$ as we have shown in the lemma. As we have seen before, if we  apply $\varepsilon$ on the second leg in the first case and on the first leg in the second case, we get
$$(\iota\ot\varepsilon)(E(1\ot a))=\varepsilon'_s(a)
\tussenen
(\varepsilon\ot\iota)((a\ot 1)E)=\varepsilon'_t(a)$$
where $\varepsilon'_s:A\to B$ and $\varepsilon'_t:A\to C$.
\snl
iii) From the proof of the lemma, we see that the range of $\varepsilon'_s$ is the same as the range of $\varepsilon_s$, namely $B$. Indeed, we have 
$$\sum_{(c)}\varepsilon_s(c_{(1)})\ot c_{(2)}aa'=\sum_{(a)} \varepsilon'_s(a_{(1)})b\ot a_{(2)}a'$$
and using the fullness of the coproduct, we see that the range of $\varepsilon_s$ is contained in the range of $\varepsilon's$. Similarly, we can define $\varepsilon'_t$ by $\varepsilon'_t(a)=(\varepsilon\ot\iota)((a\ot 1)E)$ and also $\varepsilon'_t$ and $\varepsilon_t$ have the same range, namely $C$.
\snl
iv) In the regular case, we get 
$$\varepsilon'_s(a)=\sum_{(a)} a_{(2)}S^{-1}(a_{(1)})
\tussenen
\varepsilon'_t(a)=\sum_{(a)} S^{-1}(a_{(2)})a_{(1)}$$
for $a\in A$. We see that then 
$$S(\varepsilon'_t(a))=\varepsilon_t(a)
\tussenen
S(\varepsilon'_t(a))\varepsilon_s(a)$$
for all $a$. 
\hfill $\square$ \einspr

It is somewhat remarkable that in general, the maps $\varepsilon'_s$ and $\varepsilon'_t$ exist and have the same range as the maps $\varepsilon_s$ and $\varepsilon_t$ respectively, while it is not expected that the inverse of $S$ exists. We make more comments on this peculiarity in Section 3, where we discuss further possible research.
\snl
These four counital maps are also considered in e.g.\ [B-G-L] but the notations are different. For the convenience of the reader, in [K-VD] is included an appendix with a dictionary. It includes the following formulas relating our notation with the ones used in [B-G-L]:
$$\align 
\varepsilon_s(a)&=\sqcap^R(a)
\qquad\qquad\qquad\qquad\qquad
\varepsilon_s'(a)=\overline\sqcap^R(a)\\
\varepsilon_t(a)&=\sqcap^L(a)
\qquad\qquad\qquad\qquad\qquad
\varepsilon_t'(a)=\overline\sqcap^L(a).
\endalign$$
\nl

\bf 3. Examples and special cases\rm
\nl
In this section we will treat some examples and special cases. The main purpose is to illustrate  results in Section 2 about the source and target algebras. However we will also use some of the examples for the illustration of the general theory of weak multiplier Hopf algebras because this has not yet been done in the earlier papers we wrote on the subject.
\nl
\it The groupoid examples \rm
\nl
For completeness we begin with a very brief review of the two basic motivating examples associated with a groupoid.
We will not give details as they can be found in our earlier papers on the subject (see [VD-W2] and [VD-W3]). On the other hand, we use these examples to illustrate some of the statements we made earlier in this paper, as well as for some other examples further in this section.

\inspr{3.1} Example \rm
i)
Consider a groupoid $G$. {\it First} there is the algebra $A$, defined as the space $K(G)$ of complex functions on $G$ with finite support and pointwise product. Recall that the coproduct $\Delta$ on $K(G)$ is defined by 
$$\Delta(f)(p,q)=
\cases		f(pq) & \text{if } pq \text{ is defined},\\				0 & \text{otherwise}.
\endcases$$
The pair $(A,\Delta)$ is a regular weak multiplier Hopf algebra (in the sense of Definitions 1.14 and 4.1 in [VD-W3]). The canonical idempotent  $E$ in $M(A\ot A)$  is given by the function on pairs $(p,q)$ in $G\times G$ that is $1$ if $pq$ is defined and $0$ if this is not the case. The antipode $S$ is defined by $(S(f))(p)=f(p^{-1})$ whenever $f\in K(G)$ and $p\in G$.
\snl
In this example, the algebra $A_s$ is the algebra of all complex functions on $G$ so that $f(p)=f(q)$ whenever $p,q\in G$ satisfy $s(p)=s(q)$. It is naturally identified with the algebra of all complex functions on the set $G_0$ of units in $G$. The {\it source map} $\varepsilon_s$ from $A$ to $A_s$ is defined by $(\varepsilon_s(f))(p)=f(p^{-1}p)$ whenever $p\in G$ and $f\in K(G)$. The image of the source map $\varepsilon_s(A)$, what we called in this paper the {\it source algebra}, is identified with the algebra of complex functions with finite support on the units. The  algebra $A_t$ consists of functions $f$ on $G$ so that $f(p)=f(q)$ if $t(p)=t(q)$ for $p,q\in G$. It is also identified with the space of all complex functions on the units. The {\it target map} $\varepsilon_t$ from $A$ to $A_t$ is defined by $(\varepsilon_t(f))(p)=f(pp^{-1})$ for all $p$ and $f\in K(G)$. The {\it target algebra}, i.e.\ the image $\varepsilon_t(A)$ of the target map, is again identified with the space of functions with finite support on the units. Recall that these two algebras are subalgebras of the multiplier algebra $M(A)$ (here the algebra of all complex functions on $G$). 
Observe also that the source and target algebras $\varepsilon_s(A)$ and $\varepsilon_t(A)$, can be strictly smaller than the algebras $A_s$ and $A_t$ respectively. This happens when the set of units is infinite. In that case, we see that $A_s$ is indeed the multiplier algebra $M(\varepsilon_s(A))$ of $\varepsilon_s(A)$ and similarly for the target map.
\snl
ii)
For the {\it second case}, we take the groupoid algebra $\Bbb C G$ of $G$. If we use $p\mapsto \lambda_p$ for the canonical embedding of $G$ in $\Bbb C G$, then if $p,q\in G$, we have $\lambda_p\lambda_q=\lambda_{pq}$  if $pq$ is defined and $0$ otherwise. The coproduct on $\Bbb C G$ is given by $\Delta(\lambda_p)=\lambda_p\ot\lambda_p$ for all $p\in G$. The idempotent $E$ is $\sum \lambda_e\ot \lambda_e$ where the sum is only taken over the units $e$ of $G$. The antipode is given by $S(\lambda_p)=\lambda_{p^{-1}}$ for all $p\in G$.
\snl
The {\it source} and {\it target maps} are given by $\varepsilon_s(\lambda_p)=\lambda_e$ where $e=s(p)$ and $\varepsilon_t(\lambda_p)=\lambda_e$ where now $e=t(p)$ for $p\in G$. Here the  {\it source} and {\it target algebras} coincide and it is the algebra of the span of elements of the form $\lambda_e$ where $e$ is a unit of $G$. Also here the source and target algebras need not be unital and so can be strictly smaller then their multiplier algebras.
\hfill$\square$\einspr

Recall that these two cases are dual to each other. The duality is given by $\langle f,\lambda_p\rangle=f(p)$ whenever $f\in K(G)$ and $p\in G$. We will give more details (about this duality) in  [VD-W5] where we treat duality for (regular) weak multiplier Hopf algebras with integrals.
\nl
\it Examples associated with separability idempotents \rm
\nl
For the next example, we start with any separability idempotent. Later, we will consider two special cases of this. The most important one will be constructed from the separability idempotent that is the canonical idempotent of a given weak multiplier Hopf algebra.  In some sense, we isolate the source and target algebras with what remains of the original coproduct. 
\snl
These examples illustrate very well the use of different properties of the source  and target algebras, obtained in the previous section. 
\snl 
Recall from [VD4.v2] that a  {\it separability idempotent} is an idempotent in the multiplier algebra $M(B\ot C)$ of the tensor product of two non-degenerate algebras $B$ and $C$ with certain properties. In particular there exist non-degenerate anti-homomorphisms $S_B:B\to M(C)$  and $S_C:C\to M(B)$ characterized by the formulas
$$E(b\ot 1)=E(1\ot S_B(b))
\tussenen
(1\ot c)E=(S_C(c)\ot 1)E$$
whenever $b\in B$ and $c\in C$. There are also the unique linear functionals $\varphi_B$ and $\varphi_C$ on $B$ and $C$ respectively, characterized by
$$(\varphi_B\ot \iota)(E(1\ot c))=c 
\tussenen
(\iota\ot \varphi_C)((b\ot 1)E)=b 
$$ 
for all $b\in B$ and $c\in C$. We refer to [VD4.v2] for details.
\snl
We now construct a weak multiplier Hopf algebra from a separability idempotent in the next proposition. 

\inspr{3.2} Proposition \rm 
Let $B$ and $C$ be non-degenerate algebras and assume that $E$ is a separability idempotent in $M(B\ot C)$. Let $P=C\ot B$. There is a coproduct $\Delta_P$ on $P$ defined by 
$$\Delta_P(c\ot b)=c\ot E\ot b$$
for $c\in C$ and $b\in B$. The pair $(P,\Delta_P)$ is a weak multiplier Hopf algebra. The counit $\varepsilon_P$ is given by $\varepsilon_P(c\ot b)= \varphi_B(S_C(c)b)$. We also have $\varepsilon_P(c\ot b)=\varphi_C(cS_B(b))$. The canonical idempotent $E_P$ of $(P,\Delta_P)$ is $1\ot E\ot 1$. The antipode $S_P$ is given by $S_P(c\ot b)=S_B(b)\ot S_C(c)$ when $b\in B$ and $c\in C$. 
The  source and target algebras are $1\ot B$ and $C\ot 1$ respectively and the source and target maps are
$$\varepsilon_s^P(c\ot b)=1\ot S_C(c)b
\tussenen
\varepsilon_t^P(c\ot b)=cS_B(b)\ot 1$$
for all $b\in B$ and $c\in C$. In these formulas, $1$ is the identity in $M(C)$ and $M(B)$ respectively.

\snl\bf Proof\rm:
We will systematically use $\iota_P$, $1_P$, etc.\ for objects related with $P$. For the objects related with the original algebras, we will use no index. 
\snl
i) The algebra $P$ is non-degenerate and idempotent because this is true for its components $B$ and $C$.
\snl
ii) Because $E\in M(B\ot C)$ we have that $\Delta_P(c\ot b)$, defined as $c\ot E\ot b$, belongs to $M(P\ot P)$. Because $E^2=E$, it is clear that $\Delta_P$ is a homomorphism. By assumption, we have that $E(1\ot C)$ and $(B\ot 1)E$ are subsets of $B\ot C$. Therefore 
$$\Delta_P(P)(1_P\ot P)\subseteq P\ot P
\tussenen
(P\ot 1_P)\Delta_P(P)\subseteq P\ot P.$$ 
The coproduct $\Delta_P$ is coassociative and $(\Delta_P\ot \iota_P)\Delta_P(c\ot b)=c\ot E\ot E\ot b$ for all $b\in B$ and $c\in C$.
This coproduct is full because $E$ is assumed to be full (as in Definition 1.1 of [VD4.v2]).
\snl
iii) Now we prove that there is a counit $\varepsilon_P$ on $(P,\Delta_P)$. 
First define $\varepsilon_P(c\ot b)=\varphi_C(cS_B(b))$.  
For all $b\in B$ and $c\in C$ we have that
$$\align
(\iota_P\ot\varepsilon_P)\Delta_P(c\ot b)
		&=(\iota_P\ot\varepsilon_P)(c\ot E\ot b)\\
		&=(\iota_P\ot \varphi_C)(c\ot E(1\ot S_B(b)))\\
		&=(\iota_P\ot \varphi_C)(c\ot E(b\ot 1))=c\ot b.
\endalign$$
On the other hand, if we define $\varepsilon'_P(c\ot b)=\varphi_B(S_C(c)b)$ we will find similarly
$$(\varepsilon'_P\ot \iota)\Delta_P(c\ot b)=c\ot b$$ 
for all $b\in B$ and $c\in C$. Then, from the general theory, we know that $\varepsilon_P$ and $\varepsilon'_P$ must be the same (see e.g.\ the argument we gave in the preliminary section of this paper). In the regular case we treat later, we will give another argument for this fact (see a remark after the proof of Proposition 3.3). This proves the existence of the counit.
\snl
iv) Take any elements $b,b'\in B$ and $c,c'\in C$. Then
$$\Delta_P(c\ot b)(1\ot 1\ot c'\ot b')=(1\ot E\ot 1)(c\ot 1\ot c'\ot bb').$$
If we replace $c'$ by elements of the form $S_B(b'')c''$, the right hand side will be
$$(1\ot E\ot 1)(c\ot b''\ot c''\ot bb').$$
Next we use that $B$ is idempotent and that the map $S_B$ is non-degenerate. Then we can conclude from this that 
$\Delta_P(P)(1_P\ot P)=E_P(P\ot P)$
with $E_P=1\ot E\ot 1$. Similarly we find $(P\ot 1_P)\Delta_P(P)=(P\ot P)E_P$ and it follows that $E_P$ is the canonical idempotent for $(P,\Delta_P)$.
\snl
It is  straightforward to verify that the legs of $E_P$ commute. Moreover
$$(\iota_P\ot \Delta_P)(E_P)=1\ot E\ot E\ot 1$$
and this is clearly $(1_P\ot E_P)(E_P\ot 1_P)$.
\snl
v) We now define $S_P(c\ot b)=S_B(b)\ot S_C(c)$ for all $b$ and $c$ and we show that all the conditions of Theorem 2.9 of [VD-W3] are fulfilled. This will complete the proof.
\snl
We consider the candidate for the generalized inverse $R_1$ of the canonical map $T_1$ using this expression for $S_P$. We get, using formally $E_{(1)}\ot E_{(2)}$ for $E$, that
$$\align R_1(c\ot b\ot c'\ot b')&=((\iota_P\ot S_P)(c\ot E\ot b))(1\ot 1\ot c'\ot b')\\
&=c\ot E_{(1)}\ot S_B(b)c'\ot S_C(E_{(2)})b'.
\endalign$$ 
That this maps $P\ot P$ to itself is a consequence of the property, obtained in Proposition 1.9 of [VD4.v2], saying that $E_{(1)}\ot S_C(E_{(2)})b'$ is in $B\ot B$. 
\snl
Using this candidate for the antipode, we can calculate the candidates for the source and target maps $\varepsilon^P_s$ and $\varepsilon^P_t$. We find
$$\align
\varepsilon^P_t(c\ot b)&=(c\ot E_{(1)})(S_B(b)\ot S_C(E_{(2)}))=cS_B(b)\ot 1\\
\varepsilon^P_s(c\ot b)&=(S_B(E_{(1)})\ot S_C(c))(E_{(2)}\ot b)=1\ot S_C(c)b
\endalign$$
for all $b\in B$ and $c\in C$. We have again used the Sweedler type notation for $E$ and that $E_{(1)}S_C(E_{(2)})=1$ and $S_B(E_{(1)})E_{(2)}=1$ (see Proposition 1.6 in [VD4.v2]).
\snl
Finally we have to show that 
$$\sum_{(a)}\varepsilon^P_t(a_{(1)})a_{(2)}=a
\tussenen
\sum_{(a)}\varepsilon^P_s(a_{(1)})S_P(a_{(2)})=S_P(a)$$
for all $a=c\ot b$. We find
$$\varepsilon^P_t(c\ot E_{(1)})(E_{(2)}\ot b)=cS_B(E_{(1)})E_{(2)}\ot b=c\ot b$$
proving the first equation. And 
$$\varepsilon^P_s(c\ot E_{(1)})S_P(E_{(2)}\ot b)=(1\ot S_C(c)E_{(1)})(S_B(b)\ot S_C(E_{2}))=S_B(b)\ot S_C(c).
$$
Finally, we have to show that $T_1R_1(p\ot p')=E_P(p\ot p'$ for all $p,p'\in P$ where $T_1$ is the canonical map $p\ot p'\mapsto \Delta_P(p)(1\ot p')$ and where $R_1$ it its generalized inverse constructed with the antipode $S_P$ as above. With $p=c\ot b$ and $p'=c'\ot b'$ we find
$$\align
T_1R_1(p\ot p')&=(c\ot E_{(1)}\ot E_{(2)}S_B(b)\ot 1)(1\ot 1\ot c'\ot b')\\
&=(c\ot E_{(1)}b\ot E_{(2)}\ot 1)(1\ot 1\ot c'\ot b')\\
&=(1\ot E\ot 1)(c\ot b\ot c'\ot b')
\endalign$$
and this is what we need because $E_P=(1\ot E\ot 1)$. In a similar way, we find $T_2R_2(p\ot p')=(p\ot p')E_P$ where $T_2$ is the other canonical map and $R_2$ its generalized inverse construct with the antipode $S_P$
\snl
This proves that the candidate for the antipode satisfies all the requirements needed for Theorem 2.9 of [VD-W3] and it completes the proof.
\hfill $\square$ \einspr

We now consider the regular case. The result is as expected.

\inspr{3.3} Proposition \rm If $E$ is a regular separability idempotent in $M(B\ot C)$, then the weak multiplier Hopf algebra $(P,\Delta_P)$, constructed in the previous proposition, is a regular weak multiplier Hopf algebra.

\snl\bf Proof\rm: 
There are different ways to prove this. We will use the original definitions of regularity in both cases.
\snl
Recall that $E$ is called regular if $\zeta E$ is a separability idempotent in $M(C\ot B)$ where as before $\zeta$ is the flip map. Assume that this is the case. We then have to show that the pair $(P,\Delta^{\text{cop}})$ is also a weak multiplier Hopf algebra. Here the algebra $P$ is $C\ot B$ as before while
$$\Delta^{\text{cop}}(c\ot b)=E_{(2)}\ot b\ot c\ot E_{(1)}$$
for $b\in B$ and $c\in C$. Define the isomorphism $\gamma:B\ot C\to P$ by $\gamma(b\ot c)=c\ot b$. Then the coproduct $\Delta^{\text{cop}}$ yields a coproduct $\Delta'$ on $B\ot C$ given by 
$$\Delta'(b\ot c)=b\ot E_{(2)} \ot E_{(1)}\ot c$$
for $b\in B$ and $c\in C$. Because $\zeta E$ is a separability idempotent in $M(C\ot B)$, it follows from the previous proposition that $(B\ot C,\Delta')$ is a weak multiplier Hopf algebra. Then this is also true for the pair $(P,\Delta^{\text{cop}})$. This completes the proof.
\hfill $\square$ \einspr

Observe the following. Given $b\in B$ and $c\in C$ we have
$$(S_C(c)\ot 1)E(b\ot 1)=(1\ot c)E(1\ot S_B(b))\tag"(3.1)"$$ 
and if we apply $\varphi_B\ot \varphi_C$ we find that $\varphi(S_C(c)b)=\varphi_C(cS_B(b))$. This illustrates the equality of the two forms of the counit in the formulation of Proposition 3.2. This argument however only seems to work for a (semi-)regular separability idempotent because only in that case we know that the elements in the equation (3.1) belong to $B\ot C$. 
\nl
In Proposition 3.2 we have associated a weak multiplier Hopf algebra to any separability idempotent. On the other hand, we know that conversely, the canonical idempotent $E$ of a weak multiplier Hopf algebra is a separability idempotent in $M(B\ot C)$ where now $B$ and $C$ are the source and target algebras. This is proven in Section 2 (Theorem 2.18). What happens when we then apply the construction of Proposition 3.2 again is explained the following proposition.

\inspr{3.4} Proposition \rm
Let $(A,\Delta)$ be a weak multiplier Hopf algebra. Consider the canonical idempotent $E$ as sitting in $M(B\ot C)$ where $B$ and $C$ are the source and target algebras. Associate a new weak multiplier Hopf algebra $(P,\Delta_P)$ as in Proposition 3.2. Define $\gamma:P \to M(A)$ by $\gamma(x\ot y)=xy$ for $x\in C$ and $y\in B$. Then $\gamma$ is a non-degenerate homomorphism. It satisfies $\Delta\circ\gamma=(\gamma\ot \gamma)\circ\Delta_P$ and $S\circ\gamma=\gamma\circ S_P$.

\snl\bf Proof\rm: 
Because the source and target algebras $B$ and $C$ are commuting subalgebras of $M(A)$, it follows that $\gamma$ is an algebra homomorphism from $P$ to $M(A)$. The image is $CB$. Because of Proposition 2.10, we have $CBA=A=ACB$ and so $\gamma$ is non-degenerate. It extends to a unital homomorphism on the multiplier algebra of $P$.
\snl
For all $y\in B$ and $x\in C$ we have
$$\Delta(\gamma(x\ot y))=\Delta(xy)=(x\ot y)E$$
while on the other hand
$$(\gamma\ot\gamma)\Delta_P(x\ot y)=(\gamma\ot\gamma)(x\ot E\ot y)=(x\ot 1)E(1\ot y).$$
These expressions are the same as the element $y$ commutes with the second leg of $E$.
\snl
For the antipode we find
$$\gamma(S_P(x\ot y))=\gamma(S_B(y)\ot S_C(x))=S(y)S(x)=S(xy)=S(\gamma(x\ot y))$$
where we have again used that the element $x$ of $C$ and the element $y$ of $B$ commute.
\hfill $\square$ \einspr

Remark that in general, the map $\gamma$ is not injective. Take e.g.\ the weak multiplier Hopf algebra constructed from a set $X$. The algebra $A$ is the algebra $K(X)$ of complex functions with finite support and $\Delta(f)(p, q)=0$ when $p$ and $q$ are different while $\Delta(f)(p,p)=f(p)$. This is a weak multiplier Hopf algebra. The canonical idempotent is the function $X\ot X$ that is $1$ on the diagonal and $0$ everywhere else. Clearly the left and the right legs are all of $K(X)$. In particular $B=C$. The map $\gamma$ is the multiplication map from $K(X\times X)$ to $K(X)$ and this is not injective.
\snl
If the algebra $A$ is unital, we can also show that $\gamma\circ \varepsilon_t^P=\varepsilon_t\circ\gamma$ and $\gamma\circ \varepsilon_s^P=\varepsilon_s\circ\gamma$. Indeed, for all $a$ in $A$ and $x,y$ in $C$ and $B$ respectively, we have by Proposition 2.7
$$\varepsilon_t(xya)=x\varepsilon_t(a)S(y).$$
If $a=1$ we get $\varepsilon_t(a)=1$ and so $\varepsilon_t(xy)=xS(y)$. This means $\gamma (\varepsilon_t^P(x\ot y))=xS_B(y)=\varepsilon_t(\gamma(x\ot y))$. Similarly for $\varepsilon_s$. If the algebra is not unital, we can not argue like this because the counital maps $\varepsilon_s$ and $\varepsilon_t$ have no obvious extensions from $A$ to the multiplier algebra $M(A)$. 
\snl
In [VD-W5], where we treat integrals and duality, we will consider this example again and show that integrals on $(P,\Delta_P)$ automatically exist and therefore that we can obtain a dual version of this example. 
\nl
\it Discrete quantum groups \rm
\nl
In what follows, we will use the term {\it discrete quantum group} for a regular multiplier Hopf algebra $(A,\Delta)$ of discrete type with a (left) cointegral $h$ satisfying the extra condition that $\varepsilon(h)=1$ (where $\varepsilon$ is the counit). This is the case when $h$ is an idempotent. Then $S(h)=h$ (where $S$ is the antipode) and $h$ is also a right cointegral.
\snl
It is shown in Proposition 3.11 of [VD4.v2] that $\Delta(h)$ is a separability idempotent in $M(A\ot A)$. So here $B$ and $C$ coincide with the original algebra $A$. The antipodal maps $S_B$ and $S_C$ are both nothing else but the antipode $S$ on $A$. The distinguished linear functionals $\varphi_B$ and $\varphi_C$ are the right and the left integrals $\psi$ and $\varphi$ on $(A,\Delta)$, normalized so that $\varphi(h)=\psi(h)=1$. 
\snl
Then as a consequence of Proposition 3.2, we get the following.

\inspr{3.5} Proposition \rm
Let $(A,\Delta)$ be a discrete quantum group and $h$ the normalized cointegral. The algebra $P$ defined as $A\ot A$ is a regular weak multiplier Hopf algebra for the coproduct $\Delta_P$ defined by $\Delta_P(a\ot b)=a\ot \Delta(h)\ot b$ with $a,b\in A$. The counit $\varepsilon_P$ is given by the linear map  $a\ot b\mapsto \varphi(aS(b))$. We also have $\varepsilon_P(a\ot b)=\psi(S(a)b)$. The canonical multiplier $E_P$ is $1\ot \Delta(h)\ot 1$. The antipode $S_P$ is given by $S_P(a\ot b)=S(b)\ot S(a)$ when $a,b\in A$.
The  source and target algebras are 
$$\varepsilon^P_s(P) = 1\ot A 
\tussenen
\varepsilon^P_t(P)=A\ot 1$$ 
and the source and target maps are given by 
$$\varepsilon_s^P(a\ot b)=1\ot S(a)b
\qquad\quad\text{and}\qquad\quad
\varepsilon_t^P(a\ot b)=aS(b)\ot 1$$
for all $a,b\in A$. Here $1$ is the identity in $M(A)$.
\hfill$\square$\einspr

Again we have integrals and we can construct the dual. This will be done in [VD-W5].
\nl
\it A quantization of the groupoid associated with a group action \rm
\nl
Let us start by considering the weak multiplier Hopf algebra associated with a groupoid in Example 3.1.i. We can apply the result of Proposition 3.3.
\snl
Denote the space of units by $X$. The source and target algebras $B$ and $C$ are identified with the algebra $K(X)$ of complex functions with finite support on $X$. Then the canonical idempotent is a separability idempotent in $C(X\times X)$, the algebra of all complex functions on $X\times X$. It is the function with value $1$ on the diagonal and $0$ on other elements.
\snl
We get for $P$ the algebra $K(X\times X)$ of all complex functions with finite support on $X\times X$. The element $E_P$ is the function of four variables $x,u,v,y$ in $X$ that is $1$ if $u=v$ and $0$ if $u\neq v$. The antipodal maps $S_B$ and $S_C$ on $B$ and $C$ are given by the identity map on the algebra $K(X)$. Therefore, the antipode $S_P$ on $K(X\ot X)$ is given by the flip map. 
In fact, the weak multiplier Hopf algebra we get in this way, is nothing else but the algebra of functions on the trivial groupoid $X\times X$ where the product of two elements $(x,u)$ and $(v,y)$ is only defined when $u=v$ and then is $(x,y)$.
\snl
It is also interesting to see what happens when we apply Proposition 3.4 in this case. We leave it as an excercise to the reader.
\snl
We see that this has very little to do anymore with the original groupoid. And of course, we end up with a special case of Proposition 3.2. For this, we just take any set $X$ and look at the above construction.
\nl
Let us now consider the groupoid that results from a group action on a set. So let $X$ be any set and assume that a group $H$ acts on $X$, say from the left. Denote the action as $h\tr x$ for $x\in X$ and $h\in H$. Then there is a groupoid $G$ associated as follows. One has 
$$G=\{(y,h,x) \mid x,y\in X \text{ and } h\in H \text{ so that } y=h\tr x \}.$$
The product of two elements $(z,k,y')$ and $(y,h,x)$ is defined if $y=y'$ and then
$$(z,k,y)(y,h,x)=(z,kh,x).$$
The set of units is $X$ and the source and target maps are given by
$$s(y,h,x)=x\quad\qquad\text{and}\qquad\quad t(y,h,x)=y.$$
The set of units is considered as a subset of $G$ by the embedding $x\to (x,e,x)$ where $e$ is the identity in $H$. 
\snl
We can construct the weak multiplier Hopf algebras, associated with this groupoid, as in Example 3.1. In the case where the group is trivial, we then get the example we just  mentioned above. If on the other hand, the space $X$ is trivial (i.e.\ it consists only of one point), then we get the multiplier Hopf algebras associated with the group $H$.
\nl
There is however another way to associate a weak multiplier Hopf algebra. It is a special case of a construction that we will consider next. 
\snl
The {\it starting point} is as in Proposition 3.2. We have a separability idempotent $E$ in the multiplier algebra $M(B\ot C)$ of the tensor product of two non-degenerate idempotent algebras $B$ and $C$. It need not be regular. Furthermore, we have a multiplier Hopf algebra $(Q,\Delta)$. Here we assume that it is regular. We will explain why we need this condition for the multiplier Hopf algebra.
\snl
We assume that $Q$ acts from the left on $C$ and from the right on $B$. The actions are denoted by $q\tr c$ and $b\tl q$ when $b\in B$, $c\in C$ and $q\in Q$. It is assumed that $B$ is a right $Q$-module algebra and that $C$ is a left $Q$-module algebra. In particular, the two actions are unital. Moreover, these data are required to satisfy
$$(E_{(1)}\tl q) \ot E_{(2)}=E_{(1)} \ot (q\tr E_{(2)}) \tag"(3.2)"$$
where we use the Sweedler type notation $E=E_{(1)}\ot E_{(2)}$ and where the equation is given a meaning by multiplying with an element $b$ of $B$ in the first factor from the left and with an element $c$ of $C$ in the second factor from the right.
\snl
The underlying algebra $P$ that we use in this example is a {\it two-sided smash product} of $Q$ with $B$ and $C$. The construction has been studied for Hopf algebras (see e.g.\ [B-P-VO]) but not yet for multiplier Hopf algebras. However, the results and the arguments are very similar to the theory of smash products as developed in [Dr-VD-Z]. Therefore, in the following proposition, we do not give all the details. We  concentrate on the correct statements and briefly indicate how things are proven. Remark that the construction only works fine in the case of a regular multiplier Hopf algebra. This is the reason why we need regularity for $(Q,\Delta)$.

\inspr{3.6} Proposition \rm As above, assume that $Q$ is a regular multiplier Hopf algebra, that $B$ is a right $Q$-module algebra and that $C$ is a left $Q$-module algebra. Then the tensor product $C\ot Q\ot B$ is an associative algebra $P$ with the product defined as
$$(c\ot q\ot b)(c'\ot q'\ot b')
=\sum_{(q)(q')}c(q_{(1)}\tr c')\ot q_{(2)}q'_{(1)}\ot (b\tl q'_{(2)})b'\tag"(3.3)"$$
where $b,b'\in B$, $c,c'\in C$ and $q,q'\in Q$.
\hfill $\square$\einspr

Remark that the actions are assumed to be unital and therefore they provide the necessary coverings in (3.3).
\snl
The proof of this result is straightforward. Also remark that the algebra $P$ is idempotent if this is the case for $B$ and $C$.
\snl
The two-sided smash product can be considered in two ways as a twisted product in the sense of [VD-VK]. First one considers the twisting of the algebras $C$ and $QB$ (where $QB$ is the ordinary smash product of $Q$ and $B$). In this case, the twist map is given by the formula
$$qb\ot c\mapsto \sum_{(q)} (q_{(1)}\tr c)\ot q_{(2)}b$$
where $b\in B$, $c\in C$ and $q\in Q$. For the second possibility, one takes the twisting of the algebras $CQ$ and $B$ (where $CQ$ is the smash product of $C$ and $Q$). Now the twist map is given by the formula
$$b\ot cq\mapsto \sum_{(q)}cq_{(1)}\ot (b\tl q_{(2)})$$
where again $b\in B$, $c\in C$ and $q\in Q$. In the two cases, one now has to verify that the twist map is compatible with the product in the two algebras (ensuring that the result is an associative algebra). One easily verifies that the two constructions give the same algebra and that the result is also the same as in the proposition above.
\snl
Just as in the case of smash products, one has obvious embeddings of $B$, $C$ and $Q$ in the multiplier algebra of $P$ and if we identify these three algebras with their images in $M(P)$, we see that $P$ is the linear span of elements $cqb$ with $b\in B$, $c\in C$ and $q\in Q$ and that we  have the commutation rules
\vskip 0.2 cm
\parindent 0.5cm
\item{}i)\ \ \ $B$ and $C$ commute,
\item{}ii)\ \ $bq=\sum_{(q)}q_{(1)}(b\tl q_{(2)})$ for all $b\in b$ and $q\in Q$,
\item{}iii) $qc=\sum_{(q)}(q_{(1)}\tr c) q_{(2)}$ for all $c\in C$ and $q\in Q$. 
\vskip 0.3 cm
\parindent 0 cm
Therefore we can view $P$ as the algebra generated by $B$, $C$ and $Q$ subject to these commutation rules.
\snl
By definition we have that the map $c\ot q\ot b\mapsto cqb$ is a linear bijection from $C\ot Q\ot B$ to $P$. However one also has various other maps that are also bijective. One can consider e.g.\ the maps
$$\align 
&b\ot q\ot c \mapsto bqc\\
&b\ot c\ot q \mapsto bcq\\
&q\ot b\ot c \mapsto qbc
\endalign$$
where always $b\in B$, $c\in C$ and $q\in Q$. This property will be used in the proof of Proposition 3.7 below.
\snl
Also remark that this construction reduces to well-known constructions in the following three special situations. If the multiplier Hopf  algebra $Q$ is trivial, then we obtain for $P$ simply the tensor product algebra $C\ot B$. If the algebra $B$ is trivial we obtain the smash product $C\# Q$, constructed with the right action of $Q$ on $C$ while if $C$ is trivial, we get the smash product $Q\# B$, for the left action of $B$ on $C$. Recall that in the original paper [Dr-V-Z], we developed the theory for left actions. The reader can also have a look at Section 1 of the expanded version of [De-VD-W] found on the arXiv where the two types of smash products are reviewed.
\snl
Then we are ready for the following example.

\inspr{3.7} Proposition \rm Assume that $B$ and $C$ are non-degenerate idempotent algebras and that $E$ is a  separability idempotent in $M(B\ot C)$. Let $Q$ be a  regular multiplier Hopf algebra and assume that $B$ is a right $Q$-module algebra and  $C$  a left $Q$-module algebra. Moreover assume the compatibility relation (3.2) as above.
\snl
Consider the two-sided smash product $P$ as given in the previous proposition. Then  $\Delta(q)$ and $E$ commute in the multiplier algebra of $P\ot P$ for all $q\in Q$ and the two-sided smash product $P$ can be equipped with a coproduct $\Delta_P$, defined by
$$\Delta_P(cqb)=(c\ot 1)\Delta(q)E(1\ot b) \tag"(3.4)"$$
whenever $b\in B$, $c\in C$ and $q\in Q$.
\snl
It makes of the pair $(P,\Delta_P)$ a weak multiplier Hopf algebra. The canonical idempotent $E_P$ is $E$, considered as sitting in $M(P\ot P)$. 
The counit $\varepsilon_P$ is given by the linear map  
$$qcb\mapsto \varepsilon(q)\varphi_C(cS_B(b)))$$ 
where $\varphi_C$ is the distinguished linear functional on $C$ satisfying $(\iota\ot\varphi_C)E=1$ and where $S_B$ is used for the antipodal from $B$ to $M(C)$ associated with the separability idempotent $E$. The counit $\varepsilon_P$ is also given by 
$$cbq\mapsto\varphi_B(S_B(c)b)\varepsilon(q)$$
where now $\varphi_B$ is the distinguished linear functional on $B$ and $S_B$ the antipodal map from $C$ to $M(B)$. Here $\varepsilon$ is the counit on $Q$. The antipode $S_P$ is given by $S_P(cqb)=S_B(b)S(q)S_C(c)$ when $b\in B$, $c\in C$ and $q\in Q$. Here $S$ is the antipode of the multiplier Hopf algebra $Q$.
\snl
The source and target algebras for $P$ are again the algebras $B$ and $C$, as sitting in $M(P)$ 
and the source and target maps are given by 
$$\varepsilon_s^P(cqb)=(S_C(c)\tl q)b
\tussenen
\varepsilon_t^P(cqb)=c(q\tr S_B(b))$$
for all $b\in B$, $c\in C$ and $q\in Q$. Observe that we use the extensions of the actions to the multiplier algebras.

\snl\bf Proof\rm:
First we remark that in the proof below, the coproduct, the counit, the antipode for the regular multiplier Hopf algebra $Q$ will be denoted as $\Delta$, $\varepsilon$ and $S$, without the subscript $Q$. For the coproduct, the counit and the antipode for the new weak multiplier Hopf algebra $P$, we will use subscripts and write $\Delta_P$, $\varepsilon_P$ and $S_P$. We will use superscripts for the counital maps and write $\varepsilon^P_s$ and $\varepsilon^P_t$. For the antipodal maps associated with $E$ we write $S_B$ and $S_C$. We also use $\varphi_B$ and $\varphi_C$ for the distinguished linear functionals on $B$ and $C$ respectively.
\snl
i) First, it is not hard to show that $E$ and $\Delta(q)$ for all $q\in Q$ are elements of $M(P\ot P)$. This is a consequence of the fact that the multiplier algebras of $B$, $C$ and $Q$ all sit in $M(P)$ and similarly for tensor products. 
\snl
ii) We now show that $E$ and $\Delta(q)$ commute in $M(P\ot P)$. Using the Sweedler notation, both for $E$ as before and for $\Delta(q)$ we get
$$\align E\Delta(q)
	&=\sum_{(q)} E_{(1)}q_{(1)}\ot  E_{(2)}q_{(2)}\\
	&=\sum_{(q)} q_{(1)}(E_{(1)}\tl q_{(2)})\ot  E_{(2)}q_{(3)}\\
	&=\sum_{(q)} q_{(1)}E_{(1)}\ot  (q_{(2)}\tr E_{(2)})q_{(3)}\\
	&=\sum_{(q)} q_{(1)}E_{(1)}\ot  q_{(2)}E_{(2)}=\Delta(q)E.
\endalign$$
In the above calculation, we first have used the commutation rule between $B$ and $Q$ (as the first leg of $E$ is in $B$), then the relation of the actions of $Q$ on $E$ as in formula (3.2) and finally the commutation rule between $C$ and $Q$ (as the second leg of $E$ is in $C$). Of course, to make things precise, we need to cover at the right places with the right elements. This can be done if we multiply from the left in the first factor with $bp$ and from the right in the second factor with $rc$, where $b\in B$, $c\in C$ and $p,r\in Q$.
\snl
Then we can define $\Delta_P$ on $P$ by the formula (3.4) in the formulation of the proposition. Using the commutation rules, the fact that $E$ is an idempotent, that it commutes with elements $\Delta(q)$ and that $\Delta$ is a coproduct on $Q$, it can be shown that $\Delta_P$ is a coproduct on $P$. It is full.
\snl
It is also clear that $E$, as sitting in $M(P\ot P)$ has to be the canonical idempotent for $\Delta_P$. % See also item iv) in this proof.
\snl
iii) We now prove that there is a counit and that it is given by the formulas in the formulation of the proposition. 
\snl
First define $\varepsilon_P$ on $P$ by $\varepsilon_P(qcb)=\varepsilon(q)\varphi_C(cS_B(b))$ for $b,c,q$ in $B,C,Q$ respectively. Observe that we use a different order of the elements $b,c,q$ in this definition. Then we get for all $b,c,q$ that
$$\align
(\iota_P\ot\varepsilon_P)\Delta_P(cqb)
	&=(\iota_P\ot\varepsilon_P)((c\ot 1)\Delta(q)E(1\ot b))\\
	&=\sum_{(q)} cq_{(1)}E_{(1)} \varepsilon_P(q_{(2)}E_{(2)}b)\\
	&=\sum_{(q)} cq_{(1)}E_{(1)} \varepsilon(q_{(2)})\varphi_C(E_{(2)}S_B(b))\\
	&=cqE_{(1)}b \varphi_C(E_{(2)})=cqb.
\endalign$$
If on the other hand, we  define $\varepsilon'_P$ on $P$ by the formula $\varepsilon'_P(cbq)=\varphi_B(S_C(c)b))\varepsilon(q)$, a similar calculation will give then that 
$$(\varepsilon'_P\ot\iota_P)\Delta_P(cqb)=cqb$$
for all $b,c,q$. 
\snl
It then follows from the general theory that $\varepsilon'_P=\varepsilon_P$ and that this is the counit. 
\snl
In the regular case, we consider after the proof of this proposition, we can give a direct argument for the equality of these two expressions for the counit, just as we did in the simpler case in Proposition 3.2 (see a remark after the proof of Proposition 3.3).
\snl
This takes care of the counit. 
\snl
iv) Let us now look at the {\it antipode} and the {\it source} and {\it target maps}. It is expected that the antipode $S_P$ must coincide with $S_B,S_C,S_Q$ on $B,C,Q$ respectively. 
\snl
It can be verified that $S_P$ defined in this way, is an anti-homomorphism from $P$ to $M(P)$. For this one has to argue that the definition is compatible with the commutation rules between the component $B,C,Q$. We will need to use this further in our calculations.
\snl
In order to use Theorem 2.9 of [VD-W3] again to prove that $(P,\Delta_P)$ is a weak multiplier Hopf algebra, we first must show that the candidates for the maps $R_1$ and $R_2$, constructed with the candidate for the antipode map $P\ot P$ to itself. We do this for $R_1$.
\snl
We have 
$$\align
R_1(cqb\ot c'q'b')
&=\sum_{(q)}cE_{(1)}q_{(1)}\ot S_P(E_{(2)}q_{(2)}b)c'q'b'\\
&=\sum_{(q)}cE_{(1)}q_{(1)}\ot S_B(b)S(q_{(2)})S_C(E_{(2)})c'q'b'
\endalign$$
for $c,c'\in C$, $b,b'\in B$ and $q,q'\in Q$. Then we first use that $E_{(1)}\ot S_C(E_{(2)})b''$  is in $B\ot B$ for all $b''\in B$ as we proved in Proposition 1.9 of [VD4.v2]. We use that also $\sum_{(q)}q_{(1)}\ot S(q_{(2)})q'$ is in $Q\ot Q$ for all $q,q'\in Q$. All the time, we have to shuffle elements of $B$, $C$ and $Q$ but this will not present problems. We finally get that $R_1(cqb\ot c'q'b')\in P\ot P$. The argument for $R_2$ is similar.
\snl
In order to prove the next conditions, we first calculate the candidates for the counital maps $\varepsilon^P_s$ and $\varepsilon^P_S$.  
For all $b,c,q$  we find
$$\align\varepsilon^P_s(cqb)
  &= \sum_{(q)} S_P(cE_{(1)} q_{(1)})E_{(2)}  q_{(2)}b \\
	&= \sum_{(q)} S_P(E_{(1)}c q_{(1)})E_{(2)}  q_{(2)}b\\
	&= \sum_{(q)} S(q_{(1)})S_C(c)S_B(E_{(1)}) E_{(2)}  q_{(2)}b\\
	&= \sum_{(q)} S(q_{(1)})S_C(c)  q_{(2)}b\\
	&= \sum_{(q)} S(q_{(1)})q_{(2)}(S_C(c)\tl q_{(3)})b\\
	&= (S_C(c)\tl q)b.
\endalign$$
In a similar way find 
$$\varepsilon_t^P(cqb)=c(q\tr S_B(b))$$
for all $b,c,q$ in $B,C,Q$ respectively. We use here the extension of an action to the multiplier algebra. If e.g.\ $q\in Q$ and $m\in M(C)$ we can define $q\tl m$ by the requirement $q\tl (mc)=\sum_{(q)} (q_{(1)}\tl m)q_{(2)}\tl c$ (see Proposition 4.7 in [Dr-VD-Z]).
\snl
Next we verify that $T_1R_1$ is given by left multiplication by $E$. For this, it is enough to verify that $E(cqb\ot 1)=(\iota\ot\varepsilon_t^P)\Delta_P(cqb)$ for all $b,c,q$. For the left hand side we have
$$\align
	E(cqb\ot 1)
	&=cE_{(1)}qb\ot  E_{(2)} \\
	&=\sum_{(q)} cq_{(1)} (E_{(1)}\tl q_{(2)})b \ot E_{(2)} \\
	&=\sum_{(q)} cq_{(1)} E_{(1)}b \ot q_{(2)}\tr E_{(2)} \\
	&=\sum_{(q)} cq_{(1)} E_{(1)} \ot q_{(2)}\tr (E_{(2)}S_B(b)) \\
	&=\sum_{(q)} cq_{(1)} E_{(1)} \ot (q_{(2)}\tr E_{(2)})(q_{(3)}\tr S_B(b)) \\
	&=\sum_{(q)} cq_{(1)} (E_{(1)}\tl q_{(2)})\ot  E_{(2)}(q_{(3)}\tr S_B(b)) \\
	&=\sum_{(q)} c E_{(1)}q_{(1)}\ot  E_{(2)}(q_{(2)}\tr S_B(b)) \\
	&=\sum_{(q)} E_{(1)}cq_{(1)}\ot  E_{(2)}(q_{(2)}\tr S_B(b)). 
\endalign$$
We find precisely $(\iota\ot\varepsilon_t^P)\Delta_P(cqb)$. In a similar way we find that $T_2R_2$ is given by right multiplication with $E$.
\snl
v) Finally, the only thing left is to show that 
$$\sum_{(p)}p_{(1)}S_P(p_{(2)})p_{(3)}=p
\quad\quad\text{and}\quad\quad
\sum_{(p)}S_P(p_{(1)})p_{(2)}S_P(p_{(3)})=S_P(p)$$
for all $p\in P$. We do this e.g.\ for the first one. We use that 
$$\sum_{(p)}p_{(1)}S_P(p_{(2)})p_{(3)}=\sum_{(p)}\varepsilon_t(p_{(1)})S_P(p_{(2)}).$$
Now, if $p=cqb$ we get using the Sweedler notation for $E$ that
$$\align
\sum_{(p)}\varepsilon^P_t(p_{(1)})S_P(p_{(2)})
&=\sum_{(q)}\varepsilon^P_t(cq_{(1)}E_{(1)})q_{(2)}E_{(2)}b\\
&=\sum_{(q)}cq_{(1)}\tr(S_B(E_{(1)}))q_{(2)}E_{(2)}b\\
&=\sum_{(q)}cqS_B(E_{(1)})E_{(2)}b=cqb.\endalign$$
The other formula is proven in a similarly way. This completes the proof.
\hfill $\square$\einspr

Of course, the result in Proposition 3.2 is a special case of the above. Just remark that we have to reformulate the formulas in Proposition 3.3 by considering the algebra $P$, defined as $C\ot B$ as the algebra generated by $B$ and $C$, subject to the commutation of elements of $B$ and elements of $C$ as in i) above. Elements in $P$ are then linear combinations of products $cb$ with $b\in B$ and $c$ in $C$. The coproduct $\Delta_P$ is now given as $\Delta_P(cb)=(c\ot 1)E(1\ot b)$ in $M(P\ot P)$. Also $P_s$ and $P_t$ are identified with $M(B)$ and $M(C)$, as sitting in $M(P)$ whereas the source and target maps are
$$\varepsilon_s^P(cb)=S_C(c)b 
\qquad\qquad\text{and}\qquad\qquad
\varepsilon_t^P(cb)=cS_B(b)$$
when $b\in B$ and $c\in C$.
\nl
Consider now the regular case. The following is again expected.

\inspr{3.8} Proposition  \rm
If $E$ is a regular separability idempotent, then $(P,\Delta_P)$ is a regular weak multiplier Hopf algebra.

\snl\bf Proof\rm:
We could give a direct argument as for the proof of Proposition 3.3. However, here we choose another, simpler way.
\snl
If $E$ is regular, we know that the antipodal maps $S_B$ and $S_C$ are anti-isomorphisms from $B$ to $C$ and from $C$ to $B$ respectively. Because $Q$ is also assumed to be a regular multiplier Hopf algebra, also its antipode $S$ is bijective from $Q$ to itself. This all implies that $S_P$ will map $P$ into itself and that it will be bijective. This is equivalent with saying that $(P,\Delta_P)$ is a regular weak multiplier Hopf algebra.
\hfill$\square$\einspr

We finish by giving another argument for the equality of the two expressions for the counit in the regular case.
\snl
For all $b,c,q$ we have, using again the Sweedler type notation for $E$, 
$$E_{(1)}b\ot c(q\tr E_{(2)}) =(E_{(1)}\tl q)b \ot cE_{(2)}.$$
This implies that 
$$E_{(1)}\ot c(q\tr (E_{(2)}S_B(b))) =((S_C(c)E_{(1)})\tl q)b \ot E_{(2)}.$$
If we apply $\varphi_B\ot\varphi_C$, we find
$$\varphi_C(c(q\tr S_B(b))) =\varphi_B((S_C(c)\tl q)b).\tag"(3.5)"$$
for all $b,c,q$. This is one equation we will use. 
\snl
If again we start with equation (3.2), apply $\varphi_B$ on the first factor and use fullness of $E$ we find that $\varphi_B(b\tl q)=\varepsilon(q)\varphi_B(b)$. Similarly, if we apply $\varphi_C$ on the second leg, we will get $\varphi_C(q\tr c)=\varepsilon(q)\varphi_C(c)$. In other words, the distinguished linear functionals $\varphi_B$ and $\varphi_C$ are invariant under the actions of $Q$.
\snl
Define $\varepsilon_P$ and $\varepsilon'_P$ as in the proof of Proposition 3.7. We will use the above results to give a proof of the equality of these counits. 
\snl
We find on the one hand
$$\align
\varepsilon_P(cqb)
&= \sum_{(q)}\varepsilon_P(q_{(2)}(S^{-1}(q_{(1)})\tr c)b)\\
&= \sum_{(q)}\varepsilon(q_{(2)})\varphi_C((S^{-1}(q_{(1)})\tr c)S_B(b))\\
&= \varphi_C((S^{-1}(q)\tr c)S_B(b))
\endalign$$ 
while on the other hand 
$$\align
\varepsilon'_P(cqb)
&=\sum_{(q)}\varepsilon'_P(c(b\tl S^{-1}(q_{(2)}))q_{(1)})\\
&=\sum_{(q)}\varepsilon(q_{(1)})\varphi_B(S_C(c)(b\tl S^{-1}(q_{(2)})))\\
&=\varphi_B(S_C(c)(b\tl S^{-1}(q))).
\endalign
$$
So, we need to show that 
$$\varphi_C((S^{-1}(q)\tr c)S_B(b))=\varphi_B(S_C(c)(b\tl S^{-1}(q)))\tag"(3.6)"$$
for all $b,c,q$.
\snl
For the left hand side of (3.6) we find
$$\align
\varphi_C((S^{-1}(q)\tr c)S_B(b))
&=\sum_{(q)}\varphi_C(q_{(2)}\tr((S^{-1}(q_{(1)})\tr c)S_B(b))\\
&=\varphi_C(c(q\tr S_B(b)).
\endalign$$
We have used that $\varphi_C$ is invariant under the action of $Q$. 
For the right hand side of (3.6) we get
$$\align
\varphi_B(S_C(c)(b\tl S^{-1}(q)))
&=\sum_{(q)}\varphi_B((S_C(c)(b\tl S^{-1}(q_{(2)}))\tl q_{(1)}\\
&=\varphi_B(S_C(c)\tl q)b).
\endalign$$
Here we have used that $\varphi_B$ is invariant under the action of $Q$.
\snl
Then the equation (3.6) follows from the equation (3.5)
\snl
Again, the argument does not seem to work if $E$ is not regular. Fortunately, we do not need it as we have obtained the equality in another way. 
\nl
We have not included examples of weak multiplier Hopf $^*$-algebras. In fact, the basic examples (Example 3.1) are weak multiplier Hopf $^*$-algebras for the obvious involutive structures. If in the example of Proposition 3.2, the algebras $B$ and $C$ are $^*$-algebras and if $E$ is self-adjoint, then the associated pair $(P,\Delta_P)$ will be a weak multiplier Hopf $^*$-algebra for the involutive structure on $B\ot C$ obtained from the ones on the factors $B$ and $C$. For a discrete quantum group (as in Proposition 3.5), we obtain a weak multiplier Hopf $^*$-algebra if the original discrete quantum group is a multiplier Hopf $^*$-algebra of discrete type. Finally, if in Proposition 3.7 we start with a self-adjoint separability idempotent and and with appropriate actions of a multiplier Hopf $^*$-algebra, again we will end up with a weak multiplier Hopf $^*$-algebra. 
\snl
All these statements are more or less straightforward and we leave the verification as an excercise to the reader.\nl
%\nl\nl

\newpage

\bf 4. Conclusions and further research \rm
\nl
In this paper, we have studied the source and target maps, as well as the source and target algebras, associated with a weak multiplier Hopf algebra. We have obtained results in the general case in Section 2. And we have payed special attention to the regular case. It is still not clear if the nicer results, obtained in the regular case, can be pushed forward to the non-regular case so that also there, better results can be shown. We expect however that this will not be easy, neither to prove these results if they are true, nor to find counter examples if they are not.
\snl
In fact, non-regular examples are not so easy to construct. Of course, there are the examples of Hopf algebras with a non-invertible antipode. But at this moment, we do not know of examples of multiplier Hopf algebras with a non-regular coproduct, that is with a coproduct $\Delta$ on a non-degenerate algebra $A$ so that elements of the form $\Delta(a)(b\ot 1)$ and $(1\ot c)\Delta(a)$ do not always are in $A\ot A$ for $a,b,c\in A$. More research here is needed.
\snl
Section 3 of the paper is devoted to examples. All of the examples we give are generalizations of known examples of finite-dimensional weak Hopf algebras. The duals of some of these examples, that we plan to include in [VD-W5] are probably not yet considered, even in the case of finite-dimensional weak Hopf algebras. Nevertheless, it would still be desirable to find more examples and in particular, examples that are not simply generalizations of known examples of weak Hopf algebras. We refer also to the modification procedure as explained in [VD5] to construct new examples of regular weak multiplier Hopf algebras.
\snl
The separability elements for non-unital algebras play an important role in Section 3. It is certainly worthwhile to carry out a more thorough study of these separable non-unital algebras and the associated separability idempotents (and to relate our approach with other approaches in the literature). This is partly done already in [VD4.v1]. A new version of this paper contains more information [VD4.v2]. However, there is still the open question of the existence of non-regular separability idempotents as posed in Section 5 of [VD4.v2].
\snl
Some of the examples suggest certain generalizations of the theory. Consider e.g.\ a multiplier Hopf algebra $(A,\Delta)$ of discrete type. Denote by $h$ a left cointegral. Either it can be normalized so that $\varepsilon(h)=1$ and hence $h^2=h$ (where $\varepsilon$ is the counit), or we have $\varepsilon(h)=0$ and then $h^2=0$. The first case is considered in Proposition 3.5 of this paper. The other case does not fit into this theory because $h$ and hence $\Delta(h)$ is not an idempotent. Still, it has most of the other properties of a separability idempotent. The two antipodal maps exist. Indeed, on one side we simply have
$$(1\ot a)\Delta(h)=(S(a)\ot 1)\Delta(h).$$
The other side is different because $h$ is not necessarily a right cointegral. However, by the uniqueness of cointegrals, there is a homomorphism $\gamma:A\to\Bbb C$ defined by $ha=\gamma(a)h$ for all $a$. Then
$$\Delta(h)(a\ot 1)=\Delta(h)(1\ot S'(a))$$
where $S'(a)=\sum_{(a)}\gamma(a_{(1)})S(a_{(2)})$. This is discussed in Section 5 of [VD4.v2].
\snl
Finally, as we mentioned already in the introduction, the material studied in this paper relates intimately with other research. On the one hand there is the study of weak multiplier bialgebras as introduced in [B-G-L]). We also have [K-VD] where a Larson-Sweedler type theorem is proven. Roughly it says that a weak multiplier bialgebra with enough integrals is a weak multiplier Hopf algebra. Here we have properties of the source and target maps and source and target algebras, proven in the context of weak multiplier bialgebras and separability idempotents. 
\snl
The other obvious link with the literature is the theory of multiplier Hopf algebroids as developed in [T-VD1]. In particular, there is the paper [T-VD2] where the relation between weak multiplier Hopf algebras and multiplier Hopf algebroids is studied. It seems interesting to observe that there are various possible reasons why a multiplier Hopf algebroid does not have an underlying weak multiplier Hopf algebra. 
\nl\nl

\bf References \rm
\nl
{[\bf A]} E.\ Abe: {\it Hopf algebras}. \rm Cambridge University Press (1977).
\snl
{[\bf B-G-L]} G.\ B\"ohm, J.\ G\'omez-Torecillas and E.\ L\'opez-Centella: {\it Weak multiplier bialgebras}. Trans.\ Amer.\ Math.\ Soc., ISSN 1088-6850 (online), ISSN 0002-9947 (print).  See also arXiv: 1306.1466 [math.QA]. 
\snl
{[\bf B-N-S]} G.\ B\"ohm, F.\ Nill \& K.\ Szlach\'anyi: {\it Weak Hopf algebras I. Integral theory and C$^*$-structure}. J.\ Algebra 221 (1999), 385-438. 
\snl
{[\bf B-S]} G.\ B\"ohm  \& K.\ Szlach\'anyi: {\it Weak Hopf algebras II. Representation theory, dimensions and the Markov trace}. J.\ Algebra 233 (2000), 156-212. 
\snl
{[\bf Br]} R.\ Brown: {\it From groups to groupoids: A brief survey}. Bull. London Math. Soc. 19 (1987), 113–134.
\snl
{[\bf B-P-V0]} D.\ Bulacu, F.\ Panaite \& F.\ Van Oystaeyen: {\it Generalized diagonal crossed products and smash products for quasi-Hopf algebras}. 
Commun. Math. Phys. 266 (2006), 355–399. See also arXiv:0506570 [math.QA]
\snl
{[\bf De-VD-W]} L.\ Delvaux, A.\ Van Daele \& S.\ Wang: {\it Bicrossproducts of multiplier Hopf algebras}. J.\ Algebra\ 343 (2011), 11-36.  See arXiv: 0903.2974 [math.QA] for an expanded version.
\snl
{[\bf Dr-VD]} B.\ Drabant \& A.\ Van Daele: {\it Pairing and quantum double of multiplier Hopf algebras}. 
Algebras and Representation Theory 4 (2001), 109-132. 
\snl 
{[\bf Dr-VD-Z]} B.\ Drabant, A.\ Van Daele \& Y.\ Zhang: {\it Actions of multiplier Hopf algebra}. Commun.\ Algebra 27 (1999), 4117-4172.
\snl
{[\bf H]} P.\ J.\ Higgins: {\it Notes on categories and groupoids}. Van Nostrand Reinhold, London (1971).
\snl
{[\bf K-VD]} B.-J.\ Kahng \& A.\ Van Daele: {\it The Larson-Sweedler theorem for weak multiplier Hopf algebras}. Preprint Canisius College Buffalo (USA) and University of Leuven (Belgium). See arXiv:1406.0299 [math.RA].
\snl
{[\bf N]} D.\ Nikshych: {\it On the structure of weak Hopf algebras}. Adv.\ Math.\ 170 (2002), 257-286.
\snl
{[\bf N-V1]} D.\ Nikshych \& L.\ Vainerman: {\it Algebraic versions of a finite dimensional quantum groupoid}. Lecture Notes in Pure and Applied Mathematics 209 (2000), 189-221. 
\snl
{[\bf N-V2]} D.\ Nikshych \& L.\ Vainerman: {\it Finite quantum groupiods and their applications}. In {\it New Directions in Hopf algebras}. MSRI Publications, Vol.\ 43 (2002), 211-262.
\snl
{[\bf P]} A.\ Paterson: {\it Groupoids, inverse semi-groups and their operator algebras}. Birkhauser, Boston (1999).
\snl
{[\bf R]} J.\ Renault: {\it A groupoid approach to C$^*$-algebras}. Lecture Notes in Mathematics 793, Springer Verlag.
\snl
{[\bf S]} M.\ Sweedler: {\it Hopf algebras}. Benjamin, New-York (1969).
\snl
{[\bf T]} T.\ Timmermann: {\it Regular multiplier Hopf algebroids II. Integration on and duality of algebraic quantum groupoids}. See arXiv: 1403.5282 [math.QA]
\snl 
{[\bf T-VD1]} T.\ Timmermann \& A.\ Van Daele: {\it Regular multiplier Hopf algebroids. Basic theory and examples }. See arXiv: 1307.0769 [math.QA] %... \it opgelet hier \rm ...
\snl
{[\bf T-VD2]} T.\ Timmermann \& A.\ Van Daele: {\it Multiplier Hopf algebroids arising from weak multiplier Hopf algebras}. Preprint University of M\"unster and University of Leuven. See arXiv: 1406.3509 [math.RA].
\snl 
{[\bf Va]} L.\ Vainerman (editor): {\it Locally compact quantum groups and groupoids}. IRMA Lectures in Mathematics and Theoretical Physics 2, Proceedings of a meeting in Strasbourg, de Gruyter (2002).
\snl
{[\bf VD1]} A.\ Van Daele: {\it Multiplier Hopf algebras}. Trans. Am. Math. Soc.  342(2) (1994), 917-932.
\snl
{[\bf VD2]} A.\ Van Daele: {\it An algebraic framework for group duality}. Adv.\ Math. 140 (1998), 323-366.
\snl
{[\bf VD3]} A.\ Van Daele: {\it Tools for working with multiplier Hopf algebras}.  ASJE (The Arabian Journal for Science and Engineering) C - Theme-Issue 33 (2008), 505--528. See also Arxiv: 0806.2089 [math.QA] 
\snl
{[\bf VD4.v1]} A.\ Van Daele: {\it Separability idempotents and multiplier algebras}. Preprint University of Leuven (2012). See arXiv: 1301.4398v1 [math.RA] 
\snl
{[\bf VD4.v2]} A.\ Van Daele: {\it Separability idempotents and multiplier algebras}. Preprint University of Leuven (2015). See arXiv: 1301.4398v2 [math.RA] 
\snl
{[\bf VD5]} A.\ Van Daele: {\it Modified weak multiplier Hopf algebras}. Preprint University of Leuven (2014). See arXiv:1407.0513 [math.RA].
\snl
{[\bf VD6]} A.\ Van Daele: {\it The Sweedler notation for (weak) multiplier Hopf algebras}. Preprint University of Leuven (2015). In preparation. 
\snl
{[\bf VD-VK]} A.\ Van Daele \& S.\ Van Keer: {\it The Yang-Baxter and Pentagon equation}. Comp.\ Math.\ 91 (1994), 201-221.
\snl
{[\bf VD-W1]} A.\ Van Daele \& S.\ Wang: {\it The Larson-Sweedler theorem for multiplier Hopf algebras}. J.\  Algebra\ 296 (2006), 75--95.
\snl 
{[\bf VD-W2]} A.\ Van Daele \& S.\ Wang: {\it Weak multiplier Hopf algebras. Preliminaries, motivation and basic examples}. Proceedings of  the conference 'Operator Algebras and Quantum Groups (Warsaw, September 2011), series Banach Center Publications, volume 98 (2012), 367--415. See also arXiv: 1210.3954 [math.RA] 
\snl 
{[\bf VD-W3]} A.\ Van Daele \& S.\ Wang: {\it Weak multiplier Hopf algebras I. The main theory}. Preprint University of Leuven and Southeast University of Nanjing (2012). To appear in Crelles Journal. Doi: 10.1515/crelle-2013-0053. See also Arxiv: 1210.4395 [math.RA].
\snl
{[\bf VD-W4.v1]} A.\ Van Daele \& S.\ Wang: {\it Weak multiplier Hopf algebras II. The source and target algebras}. Preprint University of Leuven and Southeast University of Nanjing (2014). For an earlier version of this paper, see arXiv:1403.7906v1 [math.RA].
\snl
{[\bf VD-W5]} A.\ Van Daele \& S.\ Wang: {\it Weak multiplier Hopf algebras III. Integrals and Duality}. University of  Leuven and Southeast University of Nanjing (in preparation). 
\end